\documentclass[12pt]{article}

\usepackage{amsmath, amsfonts, amsthm, amssymb, epsfig}

\newcommand{\R}{\mathbb{R}}
\newcommand{\N}{\mathbb{N}}
\newcommand{\Z}{\mathbb{Z}}
\newcommand{\ra}{\rightarrow}
\newcommand{\Ra}{\Rightarrow}
\newcommand{\Lra}{\Leftrightarrow}
\newcommand{\bmn}{\boldmath $n$}

\begin{document}
 \title{Exploring Progressions: A Collection of Problems}
 \author{Konstantine Zelator\\
 Department of Mathematics\\
 and Computer Science\\
 Rhode Island College \\
 600 Mount Pleasant Avenue \\
 Providence, RI  02908\\
 USA}

 \maketitle

\section{Introduction}

In this work, we study the subject of arithmetic, geometric,
mixed, and harmonic progressions. Some of the material found in
Sections 2,3,4, and 5, can be found in standard precalculus texts.
For example, refer to the books in \cite{1} and \cite{2}. A
substantial portion of the material in those sections cannot be
found in such books. In Section 6, we present 21 problems, with
detailed solutions. These are interesting, unusual problems not
commonly found in mathematics texts, and most of them are quite
challenging.  The material of this paper is aimed at mathematics
educators as well as math specialists with a keen interest in
progressions.

\section{Progressions}

In this paper we will study arithmetic and geometric progressions,
as well as mixed progressions.  All three kinds of progressions
are examples of sequences.  Almost every student who has studied
mathematics, at least through a first calculus course, has come
across the concept of sequences. Such a student has usually seen
some examples of sequences so the reader of this book has quite
likely at least some informal understanding of what the term
sequence means.  We start with a formal definition of the term
sequence.

\vspace{.15in}

\noindent{\bf Definition 1:}

\begin{enumerate}
\item[(a)] A {\bf finite sequence} of $k$ elements, ($k$ a fixed
positive integer) and  whose terms are real numbers, is a mapping
$f$ from the set $\{1,2,\ldots,k\} $ (the set containing the first
$k$ positive integers) to the set of real numbers $\R$.  Such a
sequence is usually denoted by $a_1,\ldots,a_n,\ldots,a_k$.  If
$n$ is a positive integer between $1$ and $k$, the  \bmn{\bf th
term} ${\boldmath a}_{\boldmath n}$, is simply the value of the function $f$ at
$n$; $a_{n}=f(n)$.

\item[(b)]  {\bf An infinite sequence} whose terms are real
numbers, is a mapping $f$ from the set of positive integers or
natural numbers to the set of real numbers $\R$, we write $F:\N
\ra \R$; $f(n)=a_n$.
\end{enumerate}

Such a sequence is usually denoted by $a_1,a_2,\ldots a_n,\ldots$
.  The term $a_n$ is called the $n$th term of the sequence and it
is simply the value of the function at $n$.

\vspace{.15in}

\noindent{\bf Remark 1:}  Unlike sets, for which the order in
which their elements do not  matter, in a sequence the order in
which the elements are listed does matter and makes a particular
sequence unique.  For example, the sequences $1,\ 8,\ 10,$ and
$8,\ 10,\ 1$ are regarded as different sequences.  In the first
case we have a function $f$ from $\{1,2,3\}$ to $\R$ defined as
follows:  $f: = \{1,2,3\} \ra \R;\ f(1)=1=a_1,\ f(2) = 8 = a_2$,
and $f(3)= 10= a_3$.  In the second case, we have a function $g:
\{1,2,3\} \ra \R;\ g(1)=b_1 =8,\ g(2) = b_2 = 10$, and $g(3) = b_3
=1$.

Only if two sequences are {\bf equal as functions}, are they regarded one and the same sequence.

\section{Arithmetic Progressions}

\noindent{\bf Definition 2:}  A sequence $a_1,a_2,\ldots,
a_n,\ldots $ with at least two terms, is called an {\bf
arithmetic} progression, if, and only if there exists a (fixed)
real number $d$ such that $a_{n+1} = a_n + d$, for every natural
number $n$, if the sequence is infinite.  If the sequence if
finite with $k$ terms, then $a_{n+1} = a_n+d$ for
$n=1,\ldots,k-1$. The real number $d$ is called the {\bf
difference} of the arithmetic progression.

\vspace{.15in}

\noindent{\bf Remark 2:}  What the above definition really says,
is that starting with the second term  $a_2$, each term of the
sequence is equal to the sum of the previous term plus the fixed
number $d$.

\vspace{.15in}

\noindent{\bf Definition 3:}  An arithmetic progression is said to
be {\bf increasing} if the real number $d$ (in Definition 2) is
positive, and {\bf decreasing} if the real number $d$ is negative,
and {\bf constant} if $d = 0$.

\vspace{.15in}

\noindent{\bf Remark 3:}  Obviously, if $d>0$, each term will be
greater than the previous term, while if $d < 0$, each term will
be smaller than the previous one.

\vspace{.15in}

\noindent{\bf Theorem 1:}  Let $a_1,a_2,\ldots,a_n,\ldots$ be an
arithmetic progression with difference  $d,m$ and $n$ any natural
numbers with $m < n$.  The following hold true:

\begin{enumerate}
\item[(i)]  $a_n = a_1 + (n-1)d$
\item[(ii)] $a_n = a_{n-m} + md$
\item[(iii)] $a_{m+1} + a_{n-m} = a_1 + a_n$
\end{enumerate}

\noindent {\bf Proof:}
\begin{enumerate}
\item[(i)] We may proceed by mathematical induction.  The
statement obviously holds for  $n=1$ since $a_1 = a_1 + (1-1)d$;
$a_1 = a_1 +0$, which is true.  Next we show that if the statement
holds for some natural number $t$, then this assumption implies
that the statement must also hold for $(t+1)$.  Indeed, if the
statement holds for $n=t$, then we have $a_t = a_1 + (t-1)d$, but
we also know that $a_{t+1} = a_t + d$, since $a_t$ and $a_{t+1}$
are successive terms of the given arithmetic progression.  Thus,
$a_t = a_1 + (t-1)d \Ra a_t+d = a_1(t-1)d+d \Ra a_t + d = a_1 +
d\cdot t \Ra a_{t+1} = a_1 + d \cdot t$; $a_{t+1} = a_1 + d \cdot
[(t+1) -1]$, which proves that the statement also holds for
$n=t+1$.  The induction process is complete. \item[(ii)]  By part
(i) we have established that $a_n = a_1+(n-1)d$, for every natural
number $n$.  So that

$$
\begin{array}{rcl}
a_n & = & a_1 + (n-1)d-md+ md;\\
\\
a_n & = & a_1 + [( n-m) - 1] d + md.
\end{array}
$$

Again, by part (i) we know that $a_{n-m} = a_1 + [(n-m)-1] d$.
Combining this with the  last equation we obtain, $a_n =
a_{n-m}+md$, and the proof is complete. \item[(iii)] By part (i)
we know that $a_{m+1} = a_1 + [(m+1)-1]d \Ra  a_{m+1} = a_1 + md$;
and by part (ii), we have already established that $a_n = a_{n-m}
+ md$.  Hence, $a_{m+1} + a_{n-m} = a_1 + md + a_{n-m} = a_1+a_n$,
and the proof is complete. \hfill $\square$
\end{enumerate}

\vspace{.15in}

\noindent{\bf Remark 4:}  Note that what Theorem 1(iii) really
says is that in an arithmetic  progression $a_1,\ldots, a_n$ with
$a_1$ being the first term and $a_n$ being the $n$th or last term;
if we pick two in between terms $a_{m+1}$ and $a_{n-m}$ which are
``equidistant" from the first and last term respectively
($a_{m+1}$ is $m$ places or spaces to the right of $a_1$ while
$a_{n-m}$ is $m$ spaces or places to the left of $a_n$), the sum
of $a_{m+1}$ and $a_{n-m}$ remains fixed: it is always equal to
$(a_1 +a_n)$, no matter what the value of $m$ is ($m$ can take
values from $1$ to $(n-1)$).  For example, if
$a_1,a_2,a_3,a_4,a_5$ is an arithmetic progression we must have
$a_1 +a_5 = a_2+a_4 = a_3 +a_3 = 2a_3$.  Note that $(a_2+a_4)$
corresponds to $m=1$, while $(a_3+a_3)$ corresponds to $m=2$, but
also $a_4 + a_2$ corresponds to $m=3$ and $a_5 +a_1$ corresponds
to $m=4$.

Likewise, if $b_1,b_2,b_3,b_4,b_5,b_6$ are the successive terms of an arithmetic progression we must have $b_1+b_6 = b_2 +b_5 = b_3 + b_4$.

The following theorem establishes two equivalent formulas for the sum of the first $n$ terms of an arithmetic progression.

\vspace{.15in}

\noindent{\bf Theorem 2:}  Let $a_1,a_2,\ldots ,a_n,\ldots ,$ be an arithmetic progression with difference $d$.

\begin{enumerate}
\item[(i)] The sum of the first (successive) $n$ terms
$a_1,\ldots,a_n$, is equal to the real number $\left(
{\displaystyle \frac{a_1+a_n}{2}}\right)\cdot n$; we write
$a_1+a_2+ \cdots +a_n = {\displaystyle \sum^n_{i=1}} a_i =
{\displaystyle \frac{n\cdot (a_1+a_n)}{2}}$. \item[(ii)]
${\displaystyle \sum^n_{i=1}} a_i = \left( {\displaystyle
\frac{a_1+[a_1 +(n-1)d]}{2}}\right) \cdot n$.
\end{enumerate}

\noindent{\bf Proof:}
\begin{enumerate}
\item[(i)] We proceed by mathematical induction.  For $n=1$ the statement is obviously true since $a_1 = {\displaystyle \frac{1\cdot (a_1+a_1)}{2}} = {\displaystyle \frac{2a_1}{2}}$ .  Assume the statement to hold true for some $n=k \geq 1$.  We will show that whenever the statement holds true for some value $k$ of $n,\ k\geq 1$, it must also hold true for $n=k+1$.  Indeed, assume $a_1 + \cdots + a_k = {\displaystyle \frac{k\cdot (a_1 + a_k)}{2}}$; add $a_{k+1}$ to both sides to obtain

$$\begin{array}{rcl}
a_1 + \cdots + a_k + a_{k+1}& = & {\displaystyle \frac{k\cdot a_1 +a_k}{2}} + a_{k+1}\\
& &  \Ra a_1 + \cdots + a_k + a_{k+1}\\
\\
&  = & {\displaystyle \frac{ka_1+ka_k + 2a_{k+1}}{2}}\end{array} \eqno{(1)}
$$
\noindent But since the given sequence is an arithmetic progression by Theorem 1(i), we must have $a_{k+1} = a_1+kd$ where $d$ is the difference.  Substituting back in equation (1) for $a_{k+1}$ we obtain,

$$
\begin{array}{rcl}
a_1 + \cdots + a_k + a_{k+1}& = &{\displaystyle \frac{ka_1+ka_k + (a_1+kd) + a_{k+1}}{2}}  \\
\\
\Ra a_1 + \cdots + a_k + a_{k+1}  & = & {\displaystyle \frac{(k+1)a_1 + k(a_k+d)+a_{k+1}}{2} }
\end{array}
\eqno{(2)}
$$

\noindent We also have $a_{k+1} = a_k+d$, since $a_k$ and $a_{k+1}$ are successive terms.  Replacing $a_k+d$ by $a_{k+1}$ in equation (2) we now have $a_1 + \cdots + a_k + a_{k+1} = {\displaystyle \frac{(k+1) a_1 + ka_{k+1} + a_{k+1}}{2}} = {\displaystyle \frac{(k+1)a_1+(k+1) a_{k+1}}{2}} = (k+1) \cdot {\displaystyle \frac{(a_1+a_{k+1})}{2}} $, and the proof is complete.  The statement also holds for $n=k+1$. \hfill $\square$
\item[(ii)] This is an immediate consequence of part (i).  Since ${\displaystyle \sum^n_{i=1}} a_i = {\displaystyle \frac{n(a_1+a_n)}{2}}$ and $a_n = a_1 + (n-1) d$ (by Theorem 1(i)) we have,

$${\displaystyle \sum^n_{i=1}} a_i = n\left( {\displaystyle \frac{a_1+[a_1+(n-1)d]}{2}}\right),$$

\noindent and we are done. \hfill $\square $
\end{enumerate}

\noindent {\bf Example 1:}
\begin{enumerate}
\item[(i)] The sequence of positive integers $1,2,3,\ldots ,n,\ldots ,$ is an infinite sequence which is an arithmetic progression with first term $a_1 = 1$, difference $d=1$, and the $n$th term $a_n=n$.  According to Theorem 2(i) the sum of the first $n$ terms can be easily found:  $1+2 + \ldots + n = {\displaystyle \frac{n\cdot (1+n)}{2}}$.
\item[(ii)] The sequence of the even positive integers $2,4,6,8,\ldots ,2n,\ldots $ has first term $a_1 = 2$, difference $d = 2$, and the $n$th term $a_n = 2n$.  According to Theorem 2(i), $2 + 4 + \cdots + 2n = {\displaystyle \frac{n \cdot (2 + 2n)}{2}} = {\displaystyle \frac{n\cdot 2 \cdot (n+1)}{2}} = n \cdot (n+1)$.
\item[(iii)] The sequence of the odd natural numbers $1,3,5,\ldots ,(2n-1), \ldots $, is an arithmetic progression with first term $a_1 =1$, difference $d = 2$, and $n$th term $a_n = 2n-1$.  According to Theorem 2(i) we have $1+3+ \cdots + (2n-1) = n\cdot \left( {\displaystyle \frac{1+(2n-1)}{2}}\right) = {\displaystyle \frac{n\cdot (2n)}{2}} = n^2$.
\item[(iv)] The sequence of all natural numbers which are multiples of $3\ : \ 3,6,9,12,$ \linebreak $\ldots ,3n, \ldots$ is an arithmetic progression with first term $a_1 = 3$, difference $d=3$ and $n$th term $a_n=3n$.  We have $3+6+ \cdots + 3n = {\displaystyle \frac{n\cdot (3+3n)}{2}} = {\displaystyle \frac{3n(n+1)}{2}}$.  Observe that this sum can also be found from (i) by observing that $3+6+ \cdots +3n = 3\cdot (1+2+\cdots + n) = {\displaystyle \frac{3\cdot n(n+1)}{2}}$.  If we had to find the sum of all natural numbers which are multiples of $3$, starting with $3$ and ending with $33$;  we know that $a_1 = 3$ and that $a_n = 33$.  We must find the value of $n$.  Indeed, $a_n = a_1 + (n-1)\cdot d$; and since $d = 3$, we have $33 = 3 + (n-1)\cdot 3 \Ra 33 = 3 \cdot [1+(n-1)]$; $11 = n$.  Thus, $3 + 6 + \cdots + 30 + 33  = {\displaystyle \frac{11\cdot (3 + 33)}{2}} = {\displaystyle \frac{11 \cdot 36}{2}} = 11 \cdot 18 = 198$.
\end{enumerate}

\noindent {\bf Example 2:} Given an arithmetic progression $a_1,\ldots ,a_m, \ldots , a_n, \ldots$, and natural numbers $m,n$ with $2 \leq m<n$, one can always find the sum $a_m+a_{m+1} + \cdots + a_{n-1} + a_n$;  that is, the sum of the $[(n-m)+1]$ terms starting with $a_m$ and ending with $a_n$.  If we know the values of $a_m$ and $a_n$ then we do not need to know the value of the difference.  Indeed, the finite sequence $a_m,a_{m+1}, \ldots,a_{n-1},a_n$ is a finite arithmetic progression with first term $a_m$, last term $a_n$, (and difference $d$); and it contains exactly $[(n-m)+1]$ terms.  According to Theorem 2(i) we must have $a_m+a_{m+1}+ \cdots + a_{n-1} + a_n = \frac{(n-m+1)\cdot [a_m+a_n]}{2}$.

If, on the other hand, we only know the values of the first term $a_1$ and difference $d$ ( and the values of $m$ and $n$), we can apply Theorem 2(ii) by observing that

\noindent$\begin{array}{rcl}
a_m+a_{m+1} + \cdots + a_{n-1}+a_n & = & \underset{\underset{n\ {\rm terms}}{{\rm sum\ of\ the\ first}}}{\left(\underbrace{a_1+a_2+ \ldots + a_n}\right)}\\
& & - \underset{\underset{(m-1)\ {\rm terms}}{{\rm sum\ of\ the\ first}}}{\left(\underbrace{a_1+ \ldots + a_{m-1}}\right)}\\
\\
{\rm by\ Th.\ 2(ii)}&  = & \left( \frac{2a_1 + (n-1)d}{2} \right) \cdot n\\
& &  - \left( \frac{2a_1 + (m-2)d}{2} \right) \cdot (m-1)\\
\\
& = & \frac{2[n-(m-1)]a_1 + [n\cdot (n-1)-\cdot (m-2)\cdot (m-1)]d}{2}\\
\\
& = & \frac{2(n-m+1)a_1 + [n(n-1) - (m-2)(m-1) ] d}{2}
\end{array}
$

\noindent {\bf Example 3:}

\begin{enumerate}
\item[(a)] Find the sum of all multiples of $7$, starting with $49$ and ending with $133$.
 Both $49$ and $133$ are terms of the infinite arithmetic progression with first term $a_1 = 7$, and difference $d=7$.  If $a_m = 49$, then $49 = a_1 +(m-1)d;\ 49 = 7 + (m-1) \cdot 7 \Ra \frac{49}{7} = m;\ m = 7$.  Likewise, if $a_n =$ then $133 = a_1 + (n-1)d;\ 133 = 7 + (n-1) 7 \Ra 19 = n$.  According to Example 2, the sum we are looking for is given by $a_7 + a_8 + \ldots + a_{18} + a_{19} = \frac{(19-7+1)(a_7 + a_{19})}{2} = \frac{13 \cdot (49+133)}{2} = \frac{13 \cdot 182}{2} = (13)\cdot (91) = 1183$.
\item[(b)] For the arithmetic progression with first term $a_1 = 11$ and difference $d = 5$, find the sum of its terms starting with $a_5$ and ending with $a_{13}$.

\noindent We are looking for the sum $a_5 + a_6 + \ldots + a_{12} + a_{13}$; in the usual notation $m = 5$ and $n = 13$.  According to Example 2, since we know the first term $a_1=11$ and the difference $d = 5$ we may use the formula we developed there:

$$
\begin{array}{rcl}
a_m + a_{m+1} + \ldots + a_{n-1} + a_n & = &  \frac{2(n-m+1)a_1 + [n(n-1) - (m-2)(m-1)]d}{2};\\
\\
a_5+a_6 + \ldots + a_{12} + a_{13} & = &  \frac{2\cdot(13-5+1)\cdot 11 + [13\cdot(13-1)- (5-2)(5-1)]5}{2} \\
\\
& = & \frac{2\cdot 9 \cdot 11 + [(13)(12) - (3)(4)]5}{2} =  \frac{198 + (156 -12)\cdot 5}{2} \\
\\
& = &  \frac{198 + 720}{2} =  \frac{918}{2} = 459
\end{array}
$$
\end{enumerate}

The following Theorem is simple in both its statement and proof but it serves as an effective tool to check whether three real numbers are successive terms of an arithmetic progression.

\vspace{.15in}

\noindent {\bf Theorem 3:}  Let $a,b,c$ be real numbers with $a < b < c$.
\begin{enumerate}
\item[(i)] The three numbers $a,b$, and $c$ are successive of an arithmetic progression if, and only if, $2b = a + c$ or equivalently $b = \frac{a+c}{2}$.
\item[(ii)] Any arithmetic progression containing $a,b,c$ as successive terms must have the same difference $d$, namely $d = b - a = c - b$
\end{enumerate}

\noindent{\bf Proof:}
\begin{enumerate}
\item[(i)] Suppose that $a,b$, and $c$ are successive terms of an arithmetic progression; then by definition we have $b = a + d$ and $c = b+d$, where $d$ is the difference.  So that $ d=b-a=c-b$;  from $b - a = c-b$ we obtain $2b = a+c$ or $b = \frac{a+c}{2}$.

Conversely, if $2b = a + c$, then $b - a = c-b$; so by setting $d = b - a = c - b$, it immediately follows that $b = a+d$ and $c = b+d$ which proves that the real numbers $a,b,c$ are successive terms of an arithmetic progression with difference $d$.
\item[(ii)] This has already been shown in part (i), namely that $d = b - a = c - b$.  Thus, any arithmetic progression containing the real numbers $a,b,c$ as successive terms must have difference $d = b - a = c - b$.
\end{enumerate}

\noindent {\bf Remark 5:}  According to Theorem 3, the middle term $b$ is the average of $a$ and $c$.  This is generalized in Theorem 4 below.  But, first we have the following definition.

\vspace{.15in}

\noindent {\bf Definition 4:}  Let $a_1,a_2,\ldots , a_n$ be a list (or sequence) of $n$ real numbers($n$ a positive integer).  The arithmetic {\bf mean or average} of the given list, is the real number $\frac{a_1 +a_2 + \ldots + a_n}{n}$.

\vspace{.15in}

\noindent {\bf Theorem 4:}  Let $m$ and $n$ be natural numbers with $m < n$.  Suppose that the real numbers $a_m,a_{m+1},\ldots , a_{n-1},a_n$ are the $(n-m+1)$ successive terms of an arithmetic progression (here, as in the usual notation, $a_k$ stands for the $k$th term of an arithmetic progression whose first term is $a_1$ and difference is $d$).

\begin{enumerate}
\item[(i)] If the natural number $(n-m+1)$ is odd, then the arithmetic mean or average of the reals $a_m,a_{m+1},\ldots , a_{n-1},a_n$ is the term $a_{(\frac{m+n}{2})}$.  In other words, $a_{(\frac{m+n}{2})} = \frac{a_m + a_{m+1} + \ldots + a_{n-1}+ a_n}{n-m+1}$.  (Note that since $(n-m+1)$ is odd, it follows that $n-m$ must be even, and thus so must be $n+m$; and hence $\frac{m+n}{2}$ must be a natural number).
\item[(ii)] If the natural number is even, then the arithmetic mean of the reals $a_m,a_{m+1},\ldots ,a_{n-1}, a_n$ must be the average of the two middle terms $a_{(\frac{n+m-1}{2})}$ and $a_{(\frac{n+m+1}{2})}$.

\noindent In other words $\frac{a_m+a_{m+1}+\ldots + a_{n-1} + a_n}{n-m+1} = \frac{a_{(\frac{n+m-1}{2})} + a_{(\frac{n+m+1}{2})}}{2}$.
\end{enumerate}

\noindent {\bf Remark 6:} To clearly see the workings of Theorem 4, let's look at two examples;  first suppose $m = 3$ and $n = 7$.  Then $n-m+1=7-3+1=5$; so if $a_3,a_4,a_5,a_6,a_7$ are successive terms of an arithmetic progression, clearly $a_5$ is the middle term.  But since the five terms are equally spaced or equidistant from one another (because each term is equal to the sum of the previous terms plus a fixed number, the difference $d$), it makes sense that $a_5$ would also turn out to be the average of the five terms.

If, on the other hand, the natural number $n-m+1$ is even; as in
the case of $m = 3$ and $n = 8$.  Then we have two middle numbers:
$a_5$ and $a_6$.

\vspace{.15in}

\noindent{\bf Proof} (of Theorem 4):
\begin{enumerate}
\item[(i)] Since $n-m+1$ is odd, it follows $n - m$ is even; and
thus $n+m$ is also even. Now, if we look at the integers
$m,m+1,\ldots , n-1,n$ we will see that since $m+n$ odd, there is
a middle number among them, namely the natural number
$\frac{m+n}{2}$.  Consequently among the terms $a_m,a_{m+1},\ldots
, a_{n-1},a_n$, the term $a_{(\frac{m+n}{2})}$ is the middle term.
Next we perform two calculations.  First we compute
$a_{(\frac{m+n}{2})}$ in terms of $m, n$ the first term $a_1$ and
the difference $d$.  According to Theorem 1(i), we have,
$$
a_{(\frac{m+n}{2})} = a_1 + \left( \frac{m+n}{2} - 1 \right) d = a_1 + \left( \frac{m+n-2}{2}\right) d.
$$
Now let us compute the sum $\frac{a_m+a_{m+1} + \ldots + a_{n-1} + a_n}{n-m+1}$. First assume $m \geq 2$; so that $2 \leq m < n$.  Observe that

$$\begin{array}{rl} & a_m+a_{m+1} + \ldots + a_{n-1} + a_n \\
\\
= & \underset{{\rm sum\ of\ the\ first}\ n\ {\rm terms}}{\left( {\underbrace{a_1+a_2 + \ldots + a_m + a_{m+1} + \ldots + a_{n-1}+a_n}}\right)}\\
\\
& -  \underset{\underset{{\rm note\ that}\ m-1 \geq 1,\ {\rm since}\ m \geq 2}{{\rm sum\ of\ the \ first}\ (m-1)
\ {\rm terms}}}{(\underbrace{a_1+\ldots +a_{m-1}})}
\end{array}
$$
Apply Theorem 2(ii), we have,
$$
a_1 + a_2 + \ldots + a_m + a_{m+1} + \ldots + a_{n-1} + a_n = \frac{n[2a_1 + (n-1)d]}{2}
$$
\noindent and
$$
a_1 + \ldots + a_{m-1} = \frac{(m-1)[2a_1 + (m-2)d]}{2}.
$$
Putting everything together we have
$$\begin{array}{rl}& a_m + a_{m+1} + \ldots + a_{n-1} + a_n\\
\\
 = & (a_1+a_2+\ldots +a_m + a_{m+1} + \ldots + a_{n-1} + a_n)\\
\\
& - (a_1 + \ldots + a_{m-1}) = \frac{n[2a_1+(n-1)d]}{2}\\
\\
& - \frac{(m-1)[2a_1 +(m-2)d]}{2}\\
\\
= &\frac{2(n-m+1)a_1 + [n(n-1)-(m-1)(m-2)]d}{2}.
\end{array}
$$
\noindent Thus,
$$\begin{array}{rcl}
&  & \frac{a_m+a_{m+1}+ \ldots + a_{n-1} +a_n}{n-m+1} \\
\\
& = & \frac{2(n-m+1)a_1 + [n(n-1)-(m-1)(m-2)]d}{2(n-m+1)}\\
\\
& = & a_1 + \frac{[n(n-1)-(m-1)(m-2)]d}{2(n-m+1)}\\
\\
&= & a_1 + \frac{[n^2-m^2 -n+3m-2]d}{2(n-m+1)}\\
\\
& = & a_1 + \frac{[(n-m)(n+m)+(n+m) - 2(n-m) -2]d}{2(n-m+1)}\\
\\
&= & a_1 + \frac{[(n-m)(n+m)+(n+m)-2(n-m+1)]d}{2(n-m+1)}\\
\\
&= & a_1 + \frac{[(n+m)(n-m+1) -2(n-m+1)]d}{2(n-m+1)}\\
\\
& = & a_1 + \frac{(n-m+1)(n+m-2)d}{2(n-m+1)} = a_1 + \frac{(n+m-2)d}{2},
\end{array}
$$

\noindent which is equal to the term $a_{(\frac{m+n}{2})}$ as we
have already shown. What about the case $m=1$? If $m=1$, then
$n-m+1=n$ and $a_m=a_1$.  In that case, we have the sum
$\frac{a_1+a_2+\ldots +a_{n-1}+a_n}{n} =$ (by Theorem 2(ii))
$\frac{n\cdot[2a_1+(n-1)d]}{2n}$; but the middle term
$a_{(\frac{m+n}{2})}$ is now $a_{(\frac{n+1}{2})}$ since $m = 1$;
but $a_{(\frac{n+1}{2})} = a_1 + (\frac{1+n-2}{2})d \Ra
a_{(\frac{n+1}{2})} = a_1 + (\frac{n-1}{2})d$; compare this answer
with what we just found right above, namely

$$\frac{n\cdot [2a_1 +(n-1)d]}{2n} = \frac{2a_1+(n-1)d}{2} = a_1 + (\frac{n-1}{2})d,$$

\noindent they are the same.  The proof is complete.

\item[(ii)] This is left as an exercise to the student.  (See Exercise 23).
\end{enumerate}

\vspace{.15in}

\noindent{\bf Definition 5:}  A sequence $a_1,a_2,\ldots , a_n,
\ldots$ (finite or infinite) is called a {\bf harmonic
progression}, if, and only if, the corresponding sequence of the
reciprocal terms:

$$
b_1 = \frac{1}{a_1},\ \ b_2 = \frac{1}{a_2} ,\ldots , b_n = \frac{1}{a_n}, \ldots ,
$$

\noindent is an arithmetic progression.

\vspace{.15in}

\noindent {\bf Example 4:}  The reader can easily verify that the following three sequences are harmonic progressions:

\begin{enumerate}
\item[(a)] $\frac{1}{1}, \frac{1}{2},\frac{1}{3},\ldots , \frac{1}{n}, \ldots$
\item[(b)] $\frac{1}{2}, \frac{1}{4}, \frac{1}{6} , \ldots , \frac{1}{2n}, \ldots$
\item[(c)] $\frac{1}{9}, \frac{1}{16}, \frac{1}{23} , \ldots , \frac{1}{7n+2} , \ldots$
\end{enumerate}

\section{Geometric Progressions}

\noindent {\bf Definition 6:}  A sequence $a_1,a_2,\ldots , a_n,
\ldots $ (finite or infinite) is called a {\bf geometric
progression}, if there exists a (fixed) real number $r$ such that
$a_{n+1} = r \cdot a_n$, for every natural number $n$ (if the
progression is finite with $k$ terms $a_1, \ldots , a_k$; with $k
\geq 2$, then $a_{n+1} = r\cdot a_n$, for all $n = 1,2,\ldots
,k-1$).  The real number $r$ is called the {\bf ratio} of the
geometric progression.  The first term of the arithmetic
progression is usually denoted by $a$, we write $a_1 = a$.

\vspace{.15in}

\noindent{\bf Theorem 5:}  Let $a = a_1,a_2,\ldots , a_n, \ldots$ be a geometric progression with first term $a$ and ratio $r$.
\begin{enumerate}
\item[(i)]  $a_n = a \cdot r^{n-1}$, for every natural number $n$.
\item[(ii)] $a_1 + \ldots + a_n = {\displaystyle \sum^n_{i=1}} a_i
= \frac{a_n \cdot r - a}{r-1} = \frac{a(r^n-1)}{r-1}$, for every
natural number $n$, if $r \neq 1$; if on the other hand $r = 1$,
then the sum of the first $n$ terms of the geometric progression
is equal to $n \cdot a$.
\end{enumerate}

\vspace{.15in}

\noindent{\bf Proof:}  \begin{enumerate} \item[(i)] By induction:
the statement is true for $n = 1$, since $a_1 = a \cdot r^{\circ}
= a$.  Assume the statement to hold true for $n = k$; for some
natural number $k$.  We will show that this assumption implies the
statement to be also true for $n = (k+1)$. Indeed, since the
statement is true for $n = k$, we have $a_k = a\cdot r^{k-1} \Ra
r\cdot a_k = r \cdot a\cdot r^{k-1} = a \cdot r^k$; but $k = (k+1
-1)$ and $r\cdot a_k = a_{k+1}$,  by the definition of a geometric
progression.  Hence, $a_{k+1} = a\cdot r^{(k+1)-1}$, and so the
statement also holds true for $n = k$. \item[(ii)]  Most students
probably have seen in precalculus the identity $r^n -1 =
(r-1)(r^{n-1}+ \ldots + 1)$ to hold true for all natural numbers
$n$ and all reals $r$.  For example, when $n = 2,\
r^2-1=(r-1)(r+1)$; when $n = 3$, $r^3 - 1 = (r-1)(r^2 +r + 1)$.
\end{enumerate}

We use induction to actually prove it.  Note that the statement $n
= 1$ simply takes the form, $r-1=r-1$  so it holds true; while for
$n = 2$ the statement becomes $r^2-1 = (r-1)(r+1)$, which is again
true.  Now assume the statement to hold for some $n = k,\ k \geq
2$ a natural number.  So we are assuming that the statement $r^k
-1 =(r-1)(r^{k-1} + \ldots + r + 1)$.  Multiply both sides by $r$:

$$\begin{array}{rl}
& r\cdot (r^k-1)  = r\cdot (r-1) \cdot (r^{k-1} + \ldots + r+1)\\
\\
 \Ra & r^{k+1} - r  = (r-1) \cdot (r^k+r^{k-1}+ \ldots + r^2 + r);\\
\\
 &  r^{k+1}-r  =  (r-1)\cdot (r^k + r^{k-1} +\ldots  r^2 + r + 1-1)\\
\\
 \Ra &  r^{k+1}-r = (r-1) \cdot (r^k + r^{k-1} + \ldots + r^2 + r + 1)\\
 &  + (r-1)\cdot (-1)\\
\\
 \Ra & r^{k+1} - r = (r-1) \cdot (r^k +r^{k-1} + \ldots + r^2 + r + 1) - r+1\\
\\
 \Ra & r^{k+1}-1 = (r-1) \cdot (r^{(k+1)-1} + r^{(k+1)-2} + \ldots + r^2 + r + 1),
\end{array}
$$

\vspace{.15in}

\noindent which proves that the statement also holds true for $n =
k + 1$.  The induction process is complete.  We have shown that
$r^n - 1 = (r-1)(r^{n-1}+r^{n-2} + \ldots + r + 1)$ holds true for
every real number $r$ and every natural $n$.  If $r \neq 1$, then
$r-1 \neq 0$, and so $\frac{r^n-1}{r-1} = r^{n-1} + r^{n-2} +
\ldots + r + 1$.  Multiply both sides by the first term $a$ we
obtain

$$\begin{array}{rcl}{\displaystyle \frac{a\cdot (r^n-1)}{r-1}} & = & ar^{n-1} + ar^{n-2} + \ldots ar+a\\
\\
& = & a_n + a_{n-1} + \ldots + a_2 + a_1.
\end{array}
$$

\noindent Since by part (i) we know that $a_i = a\cdot r^{i-1}$,
for $i = 1,2,\ldots , n$;  if on the other hand $r = 1$, then the
geometric progression is obviously the constant sequence,  $a,a,
\ldots ,a, \ldots\ ;\ \ a_n=a$ for every natural number $n$.  In
that case $a_1 + \ldots + a_n = \underset{n\ {\rm
times}}{\underbrace{a+ \ldots + a}} = na$.  The proof is complete.
\hfill $\square$

\vspace{.15in}

\noindent{\bf Remark 7:} We make some observation about the
different types of geometric progressions that might occur
according to the different types of values of the ratio $r$.

\begin{enumerate}
\item[(i)]  If $a = 0$, then regardless of the value of the ratio
$r$, one obtains the zero sequence $0,0,0,\ldots , 0,\ldots$ .
\item[(ii)] If $r = 1$, then for any choice of the first term $a$,
the geometric progression is the constant sequence, $a,a,\ldots ,
a,\ldots$ . \item[(iii)]  If the first term $a$ is positive and $r
> 1$ one obtains a geometric progression of positive terms, and
which is increasing and which eventually exceed any real number
(as we will see in Theorem 8, given a positive real number $M$,
there is a term $a_n$ that exceeds M; in the language of calculus,
we say that it approaches positive infinity).  For example: $a =
\frac{1}{2}$, and $r = 2$; we have the geometric progression
$$
a_1  =  a= \frac{1}{2},\ a_2 = \frac{1}{2} \cdot 2 = 1, a_3  =  \frac{1}{2} \cdot 2^2 = 2;
$$

\noindent The sequence is, $\frac{1}{2}, 1, 2, 2^2, 2^3, 2^4, \ldots , \underset{a_n}{\underbrace{\frac{1}{2} \cdot 2^{n-1}}}$.

\item[(iv)]  When $a > 0$ and $0 < r < 1$, the geometric
progression is decreasing  and in the language of calculus, it
approaches zero (it has limit value zero).

\noindent For example:  $a = 4,\ r = \frac{1}{3}$.

\noindent We have $a_1 = a = 4,\ a_2 = 4 \cdot \frac{1}{3},\ a_3 =
a\cdot  \left( \frac{1}{3}\right)^2,\ a_4 = 4 \cdot \left(
\frac{1}{3}\right)^3;$\linebreak $ 4,\ \frac{4}{3},\ \frac{4}{9} ,
\ldots , \frac{4}{3^{n-1}}\ n$th term, $ \ldots$ .

\item[(v)] For $a > 0$ and $-1 < r < 0$, the geometric sequence
alternates (which means that if we pick any term, the succeeding
term will have opposite sign).  Still, in this case, such a
sequence approaches zero (has limit value zero).

\noindent For example:  $a = 9,\ r = - \frac{1}{2}$.

\noindent $a_1 = a = 9,\ a_2 = 9\cdot \left( -\frac{1}{2}\right) =
- \frac{9}{2},\ a_3 = 9 \cdot \left(\cdot \frac{1}{2}\right)^2 =
\frac{9}{4},\ldots $

\noindent $9,\ -\frac{9}{2},\ \frac{9}{2^2}, \ \frac{-9}{2^3} ,
\ldots , \underset{n{\rm th\ term}}{\underbrace{9\cdot \left(
-\frac{1}{2} \right)^{n-1} = \frac{9\cdot (-1)^{n-1}}{2^{n-1}}}}\\
$

\item[(vi)] For $a > 0 $ and $r = -1$, we have a geometric
progression that oscillates: $a, -a,a,-a, \ldots , a_n=(-1)^{n-1},
\ldots $\ .

\item[(vii)] For $a > 0$ and $r < -1$, the geometric progression
has negative terms only, it is decreasing, and in the language of
calculus we say that approaches negative infinity.

\noindent For example:  $a = 3, r = -2$

$$\begin{array}{rcl}a_1 &  =&  a = 3,\ a_2 = 3\cdot (-2)= -6,\\
\\
 a_3 &= &3\cdot (-2)^2 = 12, \ldots3, \ -6, \ 12, \ldots ,\\
\\
&& \displaystyle{\underbrace{3 \cdot (-2)^{n-1} = 3 \cdot 2^{n-1} \cdot (-1)^{n-1}}_{n{\rm th\ term}}} , \ldots\end{array}
$$

\item[(viii)] What happens when the first term $a$ is negative?  A similar analysis holds (see Exercise 24).
\end{enumerate}

\noindent {\bf Theorem 6:}  Let $a = a_1, a_2,\ldots , a_n,\ldots $ be a geometric progression with ratio $r$.

\begin{enumerate}
\item[(i)] If $m$ and $n$ are any natural numbers such that $m <
n,\ a_n = a_{n-m} \cdot r^m$. \item[(ii)] If $m$ and $n$ are any
natural numbers such that $m < n$, then $a_{m+1}\cdot a_{n-m} =
a_1 \cdot a_n$. \item[(iii)] For any natural number $n$, $\left(
\overset{n}{\underset{i=1}{\Pi}} a_i\right) ^2 = (a_1 \cdot a_2
\ldots a_n)^2 = (a_1 \cdot a_n)^n$, where $
\overset{n}{\underset{i=1}{\Pi}} a_i$ denotes the product of the
first $n$ terms $a_1,a_2,\ldots , a_n$.
\end{enumerate}

\noindent{\bf Proof:}
\begin{enumerate}
\item[(i)] By Theorem 5(i) we have $a_n=a\cdot r^{n-1}$ and
$a_{n-m} = a\cdot r^{n-m-1}$; thus $a_{n-m} \cdot r^m = a\cdot
r^{n-m-1} \cdot r^m = a\cdot r^{n-1} = a_n$, and we are done.
\hfill $\square$ \item[(ii)] Again by Theorem 5(i) we have,

$$a_{m+1} = a \cdot r^m,\ a_{n-m} = a\cdot r^{n-m-1},\ {\rm and}\ a_n = a\cdot r^{n-1}
$$

\noindent so that $a_{m+1} \cdot a_{n-m} = a\cdot r^m \cdot a
\cdot r^{n-m-1} = a^2 \cdot r^{n-1}$ and $a_1 \cdot a_n = a\cdot
(a\cdot r^{n-1}) = a^2 \cdot r^{n-1}$.  Therefore, $a_{m+1} \cdot
a_{n-m} = a_1 \cdot a_n$. \item[(iii)] We could prove this part by
using mathematical induction.  Instead, an alternative proof can
be offered by making use of the fact that the sum of the first $n$
natural integers is equal to $\frac{n\cdot(n+1)}{2}$; $1 + 2 +
\ldots + n = \frac{n(n+1)}{2}$; we have already seen this in
Example 1(i). (Go back and review this example if necessary;
$1,2,\ldots , n$ are the consecutive first $n$ terms of the
infinite arithmetic progression with first term $1$ and difference
$1$).  This fact can be applied neatly here:

$$
\begin{array}{rcl}
a_1 \cdot a_2 \ldots a_i \ldots a_n & = & \ {\rm (by\ Theorem\ 1(i))} \\
\\
& =&  a\cdot(a\cdot r) \ldots (a\cdot r^{i-1}) \ldots (a\cdot r^{n-1})\\
\\
& = &\underset{n\ {\rm times}}{(\underbrace{a\cdot a \ldots a})}\cdot r ^{[1+2+\ldots + (i-1) +\ldots + (n-1)]}
\end{array}
$$

\noindent The sum $[1+2+\ldots + (i-1) + \ldots + (n-1)]$ is
simply the sum of the first $(n-1)$ natural numbers, if $n \geq
2$. According to Example 1(i) we have,

$$
1+2+\ldots + (i-1) + \ldots + (n-1) = \frac{(n-1)\cdot [(n-1)+1]}{2} = \frac{(n-1)\cdot n}{2}.
$$

\noindent Hence, $a_1 \cdot a_2 \ldots a_i \ldots a_n =
\underset{n\ {\rm times}}{(\underbrace{a\cdot a \ldots a})} \cdot
r^{[1+2+\ldots + (i-1) + \ldots +(n-1)]} = a^n \cdot
r^{\frac{(n-1)n}{2}} \Ra (a_1 \cdot a_2 \ldots a_i \ldots a_n)^2 =
a^{2n} \cdot r^{(n-1)\cdot n}$.  On the other hand, $(a_1 \cdot
a_n)^n = [a\cdot (a\cdot r^{n-1})]^n = [a^2 \cdot r^{n-1}]^n =
a^{2n} \cdot r^{n(n-1)} = (a_1 \cdot a_2 \ldots a_i \ldots
a_n)^2$; we are done. \hfill $\square$
\end{enumerate}

\noindent {\bf Definition 7:}  Let $a_1,a_2, \ldots, a_n$ be
positive real numbers.  The positive real number
$\sqrt[n]{a_1a_2\ldots a_n}$ is called the {\bf geometric mean} of
the numbers $a_1,a_2,\ldots a_n$.

We saw in Theorem 3 that if three real numbers $a,b,c$ are
consecutive terms of an arithmetic progression, the middle term
$b$ must be equal to the arithmetic mean of $a$ and $c$.  The same
is true for the geometric mean if the positive reals $a,b,c$ are
consecutive terms in a geometric progression.  We have the
following theorem.

\vspace{.15in}

\noindent {\bf Theorem 7:}  If the positive real numbers $a,b,c$,
are consecutive terms of a geometric progression, then the
geometric mean of $a$ and $c$ must equal $b$.  Also, any geometric
progression containing $a,b,c$ as consecutive terms, must have the
same ratio $r$, namely $r = \frac{b}{a} = \frac{c}{b}$.  Moreover,
the condition $b^2 = ac$ is the necessary and sufficient condition
for the three reals $a,b,c$ to be consecutive terms in a geometric
progression.

\vspace{.15in}

\noindent{\bf Proof:}  If $a,b,c$ are consecutive terms in a
geometric progression, then $b = ar$ and $c=b \cdot r$; and since
both $a$ and $b$ are positive and thus nonzero, we must have $r =
\frac{b}{a}={c}{b} \Ra b^2 = ac \Ra b = \sqrt{ac}$ which proves
that $b$ is the geometric mean of $a$ and $c$.  Conversely, if the
condition $b^2 = ac$ is satisfied (which is equivalent to $b =
\sqrt{ac}$, since $b$ is positive), then since $a$ and $b$ are
positive and thus nonzero, infer that $\frac{b}{a} = \frac{c}{b}$;
thus if we set $r = \frac{b}{a} = \frac{c}{b}$, it is now clear
that $a,b,c$ are consecutive terms of a geometric progression
whose ratio is uniquely determined in terms of the given reals
$a,b,c$ and any other geometric progression containing $a,b,c$ as
consecutive terms must have the same ratio $r$. \hfill $\square$

\vspace{.15in}

\noindent\framebox{
\parbox[t]{5.0in}{For the theorem to follow we will need what is called Bernoulli's Inequality:
for every real number $a \geq -1$, and every  natural number $n$,

\vspace{.15in}

\hspace{1.5in}$(a+1)^n \geq 1 + na.$}}

\vspace{.15in}

\noindent Let $a \geq -1$; Bernoulli's Inequality can be easily
proved by induction: clearly the statement holds true for $n=1$
since $1 + a \geq 1 + a$ (the equal sign holds).  Assume the
statement to hold true for some $n= k \geq 1 : (a+1)^k \geq 1 +
ka$; since $a + 1 \geq 0$ we can multiply both sides of this
inequality by $a + 1$ without affecting its orientation:

$$
\begin{array}{rcl}
(a+1)^{k+1} & \geq & (a+1)(1+ka) \Ra \\
\\
(a+1)^{k+1} & \geq & a + ka^2 +1 +ka;\\
\\
(a+1)^{k+1} & \geq & 1+ (k+1)a + ka^2 \geq 1 + (k+1)a,
\end{array}
$$

\noindent since $ka^2 \geq 0$ (because $a^2 \geq 0$ and $k$ is a natural number).  The induction process is complete.

\vspace{.15in}

\noindent{\bf Theorem 8:}
\begin{enumerate}
\item[(i)] If $r > 1$ and $M$ is a real number, then there exists
a natural number $N$ such that $r^n > M$, for every natural number
$n$. For parts (ii), (iii), (iv) and (v), let $a_1 = a, a_2,
\ldots , a_n, \ldots ,$ be an infinite geometric progression with
first term $a$ and ratio $r$.

\item[(ii)] Suppose $r > 1$ and $a > 0$.  If $M$ is a real number,
then there is a natural number $N$ such that $a_n > M$, for every
natural number $n \geq N$.

\item[(iii)] Suppose $r > 1$ and $a < 0$.  If $M$ is a real
number, then there is a natural number $N$ such that $a_n < M$,
for every natural number $n \geq N$.

\item[(iv)] Suppose $|r| < 1$, and $r \neq 0$.  If $\epsilon > 0$
is a positive real number, then there is a natural number $N$ such
that $|a_n| < \epsilon$, for every natural number $n \geq N$.

\item[(v)] Suppose $|r|<1$ and let $S_n = a_1 + a_2 + \ldots +
a_n$.  If $\epsilon > 0$ is a positive real number, then there
exists a natural number $N$ such that $\left| S_n - \frac{a}{1-r}
\right| < \epsilon$, for every natural number $n \geq N$.
\end{enumerate}

\noindent {\bf Proof:}
\begin{enumerate}
\item[(i)] We can write $r = (r-1)+1$; let $a = r-1$, since $r >
1$, $a$ must be a positive real. According to the Bernoulli
Inequality we
 have,   $r^n=(a+1)^n \geq 1+na$; thus, in order to ensure that $r^n
> M$, it is sufficient to have $1+na>M \Lra na > M-1 \Lra n > {\displaystyle \frac{M-1}{a}}$
(the last step is justified since $a > 0$).  Now, if
$\left[\hspace*{-.04in}\left[\,{\displaystyle \frac{|M-1|}{a}}
\,\right]\hspace*{-.04in}\right]$ stands for the integer part of
the positive real number ${\displaystyle \frac{|M-1|}{a}}$ we have
by definition, $\left[\hspace*{-.045in}\left[ {\displaystyle
\frac{|M-1|}{a}}\ \right]\hspace*{-.04in}\right]\ \leq \
{\displaystyle \frac{| M-1|}{a}} <
\left[\hspace*{-.04in}\left[\,{\displaystyle
\frac{|M-1|}{a}}\,\right]\hspace*{-.04in}\right]+ 1$. Thus, if we
choose $N = \left[\hspace*{-.04in}\left[ {\displaystyle
\frac{|M-1|}{a}}\right]\hspace*{-.045in}\right] +1$, it is clear
that $N > {\displaystyle \frac{|M-1|}{a}} \geq {\displaystyle
\frac{M-1}{a}}$ so that for every natural number $n \geq N$, we
will have $n > {\displaystyle \frac{M-1}{a}}$, and subsequently we
will have (since $a > 0$), $na > M-1 \Ra na + 1 > M$.  But
$(1+a)^n \geq 1 + na$ (Bernoulli), so that $r^n = (1+a)^n \geq 1 +
na > M$; $r^n > M$, for every $n \geq N$.  We are done. \hfill
$\square$ \item[(ii)] By part (i), there exists a natural number
$N$ such that $r^n > \frac{M}{a} \cdot r$, for every natural
number $n \geq N$ (apply part (i) with $\frac{M}{a}\cdot r$
replacing $M$).  Since both $r$ and $a$ are positive, so is
$\frac{a}{r}$; multiplying both sides of the above inequality by
$\frac{a}{r}$ we obtain $\frac{a}{r} \cdot r^n > \frac{a}{r} \cdot
\frac{M}{a} \cdot r \Ra a \cdot r ^{n-1} > M$.  But $a \cdot
r^{n-1}$ is the $n$th term $a_n$ of the geometric progression.
Hence $a_n > M$, for every natural number $n \geq N$. \hfill
$\square$ \item[(iii)] Apply part (ii) to the opposite geometric
progression: $-a_1,-a_2, \ldots$,\linebreak $ -a_n , \ldots $\ ,
where $a_n$ is the $n$th term of the original geometric
progression (that has  $a_1 = a < 0$ and $r>1$, it is also easy to
see that the opposite sequence is itself a geometric progression
with the same ratio $r > 1$ and opposite first term $-a$).
According to part (ii) there exists a natural number $N$ such that
$-a_n > -M$, for every natural number $n\geq N$.  Thus $-(-a_n) <
-(-M) \Ra a_n < M$, for every $n \geq N$. \hfill $\square$
\item[(iv)] Since $|r| < 1$, assuming $r \neq 0$ it follows that
$\frac{1}{|r|} > 1$.  Let \linebreak $M =
\frac{|a|}{\epsilon\cdot|r|}$.  According to part (i), there
exists a natural number $N$ such that $\left(
\frac{1}{|r|}\right)^n > M = \frac{|a|}{\epsilon |r|}$ (just apply
part (i) with $r$ replaced by $\frac{1}{|r|}$ and $M$ replaced by
$\frac{|a|}{\epsilon \cdot |r|}$ for every natural number $n \geq
N$.  Thus $\frac{1}{|r|^n} > \frac{|a|}{\epsilon \cdot |r|}$;
multiply both sides by $|r|^n \cdot \epsilon$ to obtain
$\frac{|r|^n \cdot \epsilon}{|r|^n} > \frac{|a|\cdot |r|^n\cdot
\epsilon}{\epsilon \cdot |r|} \Ra |a|\cdot |r|^{n-1} < \epsilon $;
but $|a|\cdot |r|^{n-1} = |ar^{-n} | = |a_n|$, the absolute value
of the $n$th term of the geometric progression; $|a_n| <
\epsilon$, for every natural number $n \geq N$.  Finally if $r =
0$, then $a_n = 0$ for $n \geq 2$, and so $|a_n| = 0 < \epsilon$
for all $n \geq 2$.  \hfill $\square$

\item[(v)] By Theorem 5(ii) we know that,

$$
S_n = a_1 + a_2 + \ldots + a_n = a + ar + \ldots + ar^{n-1} = \frac{a(r^n-1)}{r-1}
$$

\noindent We have $S_n - \frac{a}{1-r} = \frac{a(r^n-1)}{r-1} +
\frac{a}{r-1} = \frac{ar^n-a+a}{r-1} = \frac{ar^n}{r-1}$.
Consequently, $\left| S_n - \frac{a}{1-r}\right| = \left|
\frac{ar^n}{r-1}\right| = |r|^n \cdot \left|\frac{a}{r-1}\right|$.
Assume $r \neq 0$.  Since $|r|<1$, we can apply the already proven
part (iv), using the positive real number $\frac{\epsilon \cdot
|r-1|}{|r|}$ in place of $\epsilon$: there exists a natural number
$N$ such that $|a_n| < \frac{\epsilon\cdot|r-1|}{|r|}$, for every
natural number $n \geq N$.  But $a_n = a \cdot r^{n-1}$ so that,

$$
|a| \cdot |r|^{n-1} < \frac{\epsilon \cdot |r-1|}{|r|} \Ra
$$

\noindent $\Ra$ (multiplying both sides by $|r|$)\ \ $|a||r|^n < \epsilon \cdot |r-1| \Ra$

\noindent $\Ra$ (dividing both sides by $|r-1|$) \ \ $\frac{|a|\ |r|^n}{|r-1|} < \epsilon$.

\noindent And since $\left| S_n - \frac{a}{1-r}\right| = |r|^n
\cdot \left| \frac{a}{r-1} \right|$ we conclude that, $\left| S_n
- \frac{a}{1-r} \right| < \epsilon$.  The proof will be complete
by considering the case $r = 0$: if $r = 0$, then $a_n = 0$, for
all $n \geq 2$.  And thus $S_n = \frac{a(r^n-1)}{r-1} =
\frac{-a}{-1} = a$, for all natural numbers $n$.  Hence,  $\left|
S_n - \frac{a}{1-r} \right| = \left| a - \frac{a}{1}\right| = |a -
a| = 0 < \epsilon$, for all natural numbers $n$. \hfill $\square$
\end{enumerate}

\noindent{\bf Remark 6:}  As the student familiar with, will
recognize, part (iv) of Theorem 8 establishes the fact that the
limit value of the sequence whose $n$th term is $a_n = a \cdot
r^{n-1}$ and under the assumption $|r|< 1$, is equal to zero.  In
the language of calculus, when $|r|< 1$, the geometric progression
approaches zero.  Also, part (v), establishes the sequence of
(partial) sums whose $n$th term is $S_n$, approaches the real
number $\frac{a}{1-r}$, under the assumption $|r| < 1$.  In the
language of calculus we say that the infinite series $a + ar+ar^2
+ \ldots + ar^{n-1} + \ldots $  converges to $\frac{a}{1-r}$.

\section{\bf Mixed Progressions}

The reader of this book who has also studied calculus, may have come across the sum,

$$
1 + 2x + 3x^2 + \ldots + (n+1)x^n.
$$

\noindent There are $(n+1)$ terms in this sum; the $i$th term is
equal to $i \cdot x^{i-1}$, where $i$ is a natural number between
$1$ and $(n+1)$.  Note that if $a_i = i \cdot x^{i-1},\ b_i = i$,
and $c_i = x^{i-1}$, we have $a_i = b_i \cdot c_i$; what is more,
$b_i$ is the $i$th term of an arithmetic progression (that has
both first term and difference equal to $1$); and $c_i$ is the
$i$th term of a geometric progression (with first term $c = 1$ and
ratio $r = x$).  Thus the term $a_i$ is the product of the $i$th
term of an arithmetic progression with the $i$th term of a
geometric progression; then we say that $a_i$ is the $i$th term of
a mixed progression.  We have the following definition.

\vspace{.15in}

\noindent{\bf Definition 8:}  Let $b_1, b_2, \ldots , b_n, \ldots$
be an arithmetic progression; and $c_1,c_2,$  \linebreak $\ldots ,
c_n, \ldots$ be a geometric progression.  The sequence $a_1, a_2,
\ldots, a_n, \ldots ,$ where $a_n = b_n \cdot c_n$, for every
natural number $n$,  is called {\bf a mixed progression}.  (Of
course, if both the arithmetic and geometric progressions are
finite sequences with the same number of terms, so it will be with
the mixed progression.)

Back to our example. With a little bit of ingenuity, we can compute this sum;  that is, find a closed form expression for it, in terms of $x$ and $n$.  Indeed, we can write the given sum in the form,
$$
\begin{array}{ll}\underset{(n+1)\ {\rm terms}}{\left(\underbrace{1 + x + x^2 + \ldots + x^{n-1} + x^n}\right)} + \underset{n\ {\rm terms}}{\left( \underbrace{x + x^2 + \ldots + x^{n-1}+x^n}\right)}\\
\\
+\underset{(n-1)\ {\rm terms}}{\left( \underbrace{x^2 + x^3 + \ldots + x^{n-1}+x^n}\right)}+ \ldots + \underset{2\ {\rm terms}}{\left(\underbrace{x^{n-1}+ x^n}\right)} + \underset{{\rm one\ term}}{\underbrace{x^n}}\end{array}.$$
\noindent In other words we have written the original sum $1 + 2x + 3x^2 + \ldots + (n+1)x^n$  as a sum of $(n+1)$ sums, each containing one term less than the previous one.

Now according to Theorem 5(ii),
$$
1+x+x^2 + \ldots + x^{n-1} + x^n = \frac{x^{n+1}-1}{x-1}\ {\rm \left(\right.}{\rm assuming}\   x\neq 1\ {\rm \left. \right)},
$$
\noindent since this is the sum of the first $(n+1)$ terms of a geometric progression with first term $1$ and ratio $x$.

Next, consider
$$\begin{array}{rcl} x + x^2 + \ldots + x^{n-1} + x^n &=& (1+x+x^2+ \ldots + n^{n-1} + x^n) -1 \\
 & = & \frac{x^{n+1}-1}{x-1} - \left( \frac{x^i -1}{x-1}\right)\end{array}.
$$
\noindent Continuing this way we have,
$$\begin{array}{rcl}
x^2 + \ldots + x^{n-1} + x^n & = & (1+x+x^2 + \ldots + x^n-1 + x^n) - (x + 1) \\
\\
& = & {\displaystyle \frac{x^{n+1}-1}{x-1}} - \left({\displaystyle \frac{x^2-1}{x-1}}\right) .
\end{array}
$$
\noindent On the $i$th level,
$$
\begin{array}{rcl}
x^i + \ldots + x^{n-1} + x^n & = & (1+x+\ldots +x^{i-1} + x^i+\ldots + x^{n-1} + x^n)-\\
\\
-(1+x+\ldots + x^{i-1}) & = & {\displaystyle \frac{x^{n+1}-1}{x-1}} - \left({\displaystyle \frac{x^i-1}{x-1}}\right) .
\end{array}
$$
Let us list all of these sums:
$$
\begin{array}{crcl}
(1) & 1 + x + x^2 + \ldots + x^{n-1} + x^n & = & \frac{x^{n+1}-1}{x-1} \\
\\
(2) & x+ x^2 + \ldots + x^{n-1} + x^n & = & \frac{x^{n+1}-1}{x-1} - \left( \frac{x-1}{x-1}\right)\\
\\
(3) & x^2 + \ldots + x^{n-1} + x^n & = & \frac{x^{n+1}-1}{x-1} - \left( \frac{x^2-1}{x-1}\right)\\
\vdots & & & \\
(i) & x^i + \ldots + x^{n-1}+x^n & = & \frac{x^{n+1}-1}{x-1} - \left( \frac{x^i -1}{x-1}\right)\\
\vdots & & & \\
(n) & x^{n-1}+x^n & = & \frac{x^{n+1}-1}{x-1} - \left( \frac{x^{n-1}-1}{x-1}\right)\\
\\
(n+1) & x^n & = & \frac{x^{n+1}-1}{x-1} - \left( \frac{x^n-1}{x-1}\right),
\end{array}
$$
\noindent with $x \neq 1$.

If we add the $(n+1)$ equations or identities (they hold true for
all reals except for $x = 1$), the sum of the $(n+1)$ left-hand
sides is simply the original sum $1+2x+3x^2 + \ldots + nx^{n-1} +
(n+1)x$?  Thus, if we add up the $(n+1)$ equations member-wise we
obtain,
\newpage

$$
\begin{array}{rl}
& 1+2x+3x^2 + \ldots + nx^{n-1} + (n+1)x^n \\
= & (n+1) \cdot \left(\frac{x^{n+1}-1}{x-1}\right) + \frac{n-(x+x^2 + \ldots + x^i + \ldots + x^{n-1} + x^n)}{x-1} \\
 = & (n+1) \cdot \left( \frac{x^{n+1}-1}{x-1}\right) + \frac{(n+1) - (1 + x + x^2 + \ldots + x^n)}{x-1}\\
 \Ra & 1 + 2x + 3x^2 + \ldots + nx^{n-1} + (n+1)x^n \\
 = & (n+1)\cdot \left(\frac{x^{n+1}-1}{x-1}\right) +\frac{(n+1)-\left(\frac{x^{n+1}-1}{x-1}\right)}{x-1};\\
 & 1 + 2x +3x^2 + \ldots + nx^{n-1} + (n+1)x^n \\
= & (n+1) \cdot \left( \frac{x^{n+1}-1}{x-1} \right) + \frac{(n+1)(x-1)-(x^{n+1}-1)}{(x-1)^2};\\
& 1 + 2x+3x^2 + \ldots + nx^{n-1}+(n+1)x^n\\
 = &\frac{(n+1)(x^{n+1}-1) (x-1) + (n+1)(x-1) - (x^{n+1}-1)}{(x-1)^2};\\
& 1 + 2x + 3x^2 + \ldots + nx^{n-1}+ (n+1)x^n\\
 = & \frac{(n+1)(x-1)\cdot [(x^{n+1}-1)+1]  - (x^{n+1}-1)}{(x-1)^2};\\
& 1 + 2x + 3x^2 + \ldots + nx^{n-1}+ (n+1)x^n \\
 = & \frac{(n+1) (x-1) \cdot x^{n+1}- (x^{n+1}-1)}{(x-1)^2};\\
& 1 + 2x + 3x^2 + \ldots + nx^{n-1}+ (n+1)x^n\\
 = & \frac{(n+1)x^{n+2}-(n+1)x^{n+1}-x^{n+1} +1 }{(x-1)^2};\\
& 1 + 2x + 3x^2 + \ldots + nx^{n-1}+ (n+1)x^n\\
 = & \fbox{$\frac{(n+1)x^{n+2}-(n+2)x^{n+1} +1}{(x-1)^2}$}
\end{array}
$$

\noindent for every natural number $n$.

For $x = 1$, the above derived formula is not valid.  However, for
$x =1 $; $1 + 2x + 3x^2 + \ldots + nx^{n-1}+(n+1)x^n = 1 + 2 + 3 +
\ldots + n + (n+1) = \frac{(n+1)(n+2)}{2}$ (the sum of the first
$(n+1)$ terms of an arithmetic progression with first term $a_1 =
1$ and difference $d = 1$.

The following theorem gives a formula for the sum of the first $n$ terms of a mixed progression.

\vspace{.15in}

\noindent{\bf Theorem 9:}  Let $b_1,b_2, \ldots ,b_n , \ldots$ ,
be an arithmetic progression with first term $b_1$ and difference
$d$;  and $c_1,c_2, \ldots , c_n, \ldots $ , be a geometric
progression with first term $c_1 = c$ and ratio $r \neq 1$.  Let
$a_1, a_2,\ldots , a_n , \ldots$ , the corresponding mixed
progression, that is the sequence whose $n$th term $a_n$ is given
by $a_n = b_n \cdot c_n$, for every natural number $n$.

\begin{enumerate}
\item[(i)]  $a_n = \left[ b_1 + (n-1) \cdot d\right] \cdot c \cdot r^{n-1}$, for every natural number $n$.
\item[(ii)] For every natural number $n$, $a_{n+1} - r\cdot a_n = d \cdot c_{n+1}$.
\item[(iii)] If $S_n = a_1 + a_2 + \ldots + a_n$ (sum of the first $n$ terms of the mixed progression), then

$$
\begin{array}{rcl}
S_n & = &  \frac{a_n \cdot r - a_1}{r-1} + \frac{d\cdot \tau \cdot c \cdot (1-r^{n-1})}{(r-1)^2};\\
\\
S_n & = & \frac{a_n \cdot r - a_1}{r-1} + \frac{d\cdot r \cdot (c - c_n)}{(r-1)^2}
\end{array}
$$

\noindent (recall $c_n = c\cdot r ^{n-1}$).
\end{enumerate}

\noindent{\bf Proof:}
\begin{enumerate}
\item[(i)] This is immediate, since by Theorem 1(i), $b_n = b_1 +
(n-1)\cdot d$ and by Theorem 5(i), $c_n = c\cdot r^{n-1}$, and so
$a_n = b_n \cdot c_n = [ b_1 + (n-1)d]\cdot c \cdot r ^{r-1}$.
\item[(ii)] We have $a_{n+1} = b_{n+1} \cdot c_{n+1},\ a_n = b_n
c_n,\ b_{n+1} = d + b_n$.  Thus, $a_{n+1} - r \cdot a_n = c_{n+1}
\cdot (d + b_n) - r\cdot b_n \cdot c_n = d\cdot c_{n+1} +
c_{n+1}b_n - rb_n c_n = d \cdot c_{n+1} + b_n \cdot
\underset{0}{(\underbrace{c_{n+1} -rc_n})}= dc_{n+1}$, since
$c_{n+1}= rc_n$ by virtue of the fact that $c_n$ and $c_{n+1}$ are
consecutive terms of a geometric progression with ratio $r$.  End
of proof. \hfill $\square$ \item[(iii)] We proceed by mathematical
induction.  The statement is true for $n = 1$ because $S_1 = a_1$
and $\frac{a_1r-a_i}{r-1} + \frac{d\cdot r \cdot (c-c_1)}{(r-1)^2}
= \frac{a_1(r-1)}{r-1} + 0 = a_1 = S_1$.  Assume the statement to
hold for $n=k$:  (for some natural number $k \geq 1; \ S_k =
\frac{a_k \cdot r - a_1}{r-1} + \frac{d\cdot r \cdot (c-
c_k)}{(r-1)^2}$.  We have $S_{k+1} = S_k + a_{k+1} = \frac{a_k
\cdot r-a_1}{r-1} + \frac{d \cdot r \cdot (c - c_k)}{(r-1)^2} +
a_{k+1} = \frac{a_k \cdot r-a_1 + a_{k+1} \cdot r - a_{k+1}}{r-1}
+ \frac{d \cdot r \cdot (c - c_k)}{(r-1)^2}  (1)$. But by part
(ii) we know that $a_{k+1} - ra_k = d\cdot c_{k+1}$.  Thus, by (1)
we now have,

$$
\begin{array}{lrcl}
& S_{k+1} & = & {\displaystyle \frac{a_{k+1}\cdot r - a_1}{r-1} - \frac{d\cdot c_{k+1}}{r-1} + \frac{d\cdot r \cdot (c - c_k)}{(r-1)^2}}\\
\Ra & S_{k+1}& =& {\displaystyle \frac{a_{k+1} \cdot r - a_1}{r-1} + \frac{-(r-1)\cdot d \cdot c_{k+1} + d \cdot r\cdot (c-c_k)}{(r-1)^2}};\\
& S_{k+1} & = & {\displaystyle \frac{a_{k+1}\cdot r-a_1}{r-1} +
\frac{d\cdot r\cdot (c - c_{k+1}) + d\cdot
\overset{0}{(\overbrace{c_{k+1}- r \cdot c_k})}}{(r-1)^2}}.
\end{array}
$$

\noindent But $c_{k+1} - r\cdot c_k = 0$ (since $c_{k+1} = r \cdot
c_k$) because $c_k$ and $c_{k+1}$ consecutive terms of a geometric
progression with ratio $r$.  Hence, we obtain $S_{k+1} =
\frac{a_{k+1}\cdot r-a_1}{r-1} + \frac{d \cdot r \cdot (c -
c_{k+1})}{(r-1)^2}$; the induction is complete.
\end{enumerate}

The example with which we opened this section is one of a mixed
progression.  We dealt with the sum $1 + 2x+3x^2 + \ldots +
nx^{n-1} + (n+1)x^n$.  This is the sum of the first $(n+1)$ terms
of a mixed progression whose $n$th term is $a_n=n\cdot x^{n-1}$;
in the notation of Theorem 9, $b_n = n,\ d=1,\ c_n = x^{n-1}$, and
$r = x$ (we assume $x \neq 1$).

According to Theorem 9(iii)

$$
\begin{array}{rcl}
S_n & = & 1 + 2x+3x^2 + \ldots + nx^{n-1} = \frac{(nx^{n-1}) \cdot x - 1}{x-1} + \frac{x\cdot(1 - x^{n-1})}{(x-1)^2} \\
\\
 & = & \frac{nx^n-1}{x-1} + \frac{x-x^n}{(x-1)^2} = \frac{(nx^n  -1)(x-1)}{(x-1)^2} + \frac{x-x^n}{(x-1)^2 }\\
\\
& = & \frac{nx^{n+1}-nx^n - x+1 + x - x^n}{(x-1)^2} = \frac{nx^{n+1} - (n+1) x^n + 1}{(x-1)^2};
\end{array}
$$

\noindent Thus, if we replace $n$ by $(n+1)$ we obtain, $S_{n+1} =
1 + 2x + 3x^2 + \ldots + nx^{n-1} +(n+1)x^n =
\frac{(n+1)x^{n+2}-(n+2)x^{n+1}+1}{(x-1)^2}$,  and this is the
formula we obtained earlier.

\vspace{.15in}

\noindent{\bf Definition 9:}  Let $a_1, \ldots , a_n$ be nonzero
real numbers.  The real number  $\frac{n}{\frac{1}{a_1} + \ldots +
\frac{1}{a_n}}$, is called the {\bf harmonic mean}  of the real
numbers $a_1, \ldots, a_n$.

\vspace{.15in}

\noindent {\bf Remark 7:}  Note that since $\frac{n}{\frac{1}{a_1}
+ \ldots + \frac{1}{a_n}} = \frac{1}{(\frac{1}{a_1} + \ldots +
\frac{1}{a_n})/n}$, the harmonic mean of the reals $a_1, \ldots ,
a_n$, is really the reciprocal of the mean of the reciprocal real
numbers $\frac{1}{a_1}, \ldots , \frac{1}{a_n}$.

\vspace{.15in}

We close this section by establishing an interesting, significant
and deep inequality, that has many applications  in mathematics
and has been used to prove a number of other theorems.  Given $n$
positive real numbers $a_1, \ldots , a_n$ one can always designate
three positive reals to the given set $\{ a_1, \ldots , a_n \}$:
the arithmetic mean denoted by A.M., the geometric mean denoted by
G.M., and the harmonic mean H.M.  The
arithmetic-geometric-harmonic mean inequality asserts that A.M.
$\geq$ G.M. $\geq$ H.M. (To the reader:  Do an experiment; pick a
set of three positive reals; then a set of four positive reals;
for each set compute the A.M., G.M., and H.M. values; you will see
that the inequality holds; if you are in disbelief do it again
with another sample of positive real numbers.)

The proof we will offer for the arithmetic-geometric-harmonic
inequality is indeed short.  To do so, we need a preliminary
result: we have already proved (in the proof of Theorem 5(i)) the
identity $r^n -1 = (r-1)(r^{n-1}+ r^{n-2}+ \ldots + r+1)$, which
holds true for all real numbers $r$ and all natural numbers $n$.
Moreover, if $r \neq 1$, we have

$$
\frac{r^{n-1}}{r-1} = r^{n-1}+ r^{n-2} + \ldots + r + 1
$$

\noindent If we set $ r = \frac{b}{a}$, with $b \neq a$, in the above equation and we multiply both sides by $a^n$ we obtain,

$$
\frac{b^n - a^n}{b-a} = b^{n-1} + b^{n-2} \cdot a + b^{n-3} \cdot a^2 + \ldots + b^2 \cdot a^{n=-3} + b \cdot a^{n-2} + a^{n-1}
$$

\noindent Now, if $b > a > 0$ and in the above equation we replace
$b$ by $a$, the resulting right-hand side will be smaller.  In
other words, in view of $b > a > 0$ we have,

\vspace{.15in}

\ \hspace*{-.25in}$\begin{array}{c}
(1)\\
(2)\\
(3)\\
\vdots\\
(n-2)\\
(n-1)\\
(n)
\end{array}
\left\{ \begin{array}{l} b^{n-1} > a^{n-1}\\
b^{n-2}\cdot a > a^{n-2}\cdot a^1 = a^{n-1}\\
b^{n-3} \cdot a^2 > a^{n-3}\cdot a^2 = a^{n-1}\\
\vdots \\
b^2 \cdot a^{n-3} \cdot a^2 > a^2 \cdot a^{n-3} j= a^{n-1}\\
b \cdot a^{n-2} > a \cdot a ^{n-2} = a^{n-1}\\
a^{n-1} = a^{n-1}\end{array}\right\} \begin{array}{ll}\Ra&{\rm add\ memberwise}\\
\\
& b^{n-1} + b^{n-2} \cdot a+b^{n-3} \cdot a^2 + \ldots \\
+&  b^2a^{n-3}+b\cdot a^{n-2} + a^{n-1}\\
>&  n\cdot a^{n-1}\end{array}$

\vspace{.15in}

\noindent Hence, the identity above, for $b > a > 0$, implies the
inequality $\frac{b^n-a^n}{b-a}>na^{n-1}$; multiplying both sides
by $b-a > 0$ we arrive at

$$
\begin{array}{rl}
& b^n-a^n > (b-a)na^{n-1}\\
\\
\Ra & b^n > nba^{n-1} - na^n + a^n;\\
\\
& b^n > nba^{n-1} - (n-1)a^n.
\end{array}
$$

Finally, by replacing $n$ by $(n+1)$ in the last inequality we obtain,

\vspace{.15in}

\fbox{$b^{n+1} > (n+1)ba^n - na^{n+1}$, \begin{tabular}{l}for every natural number $n$ and any\\
real numbers such that $b>a>0$\end{tabular}}

\vspace{.15in}

We are now ready to prove the last theorem of this chapter.

\vspace{.15in}

\noindent {\bf Theorem 10:}  Let $n$ be a natural number and $a_1 , \ldots , a_n$ positive real numbers.  Then,

$$\begin{array}{rcccl}
\underset{{\rm A.M.}}{\underbrace{\frac{a_1\ldots + a_n}{n}}} & \geq & \underset{{\rm G.M.}}{\underbrace{\sqrt[n]{a_1\ldots a_n}}} & \geq &
\underset{{\rm H.M.}}{\underbrace{\frac{n}{\frac{1}{a_1} +\frac{1}{a_2}+\ldots +\frac{1}{a_n}}}}
\end{array}
$$

\noindent{\bf Proof:}  Before we proceed with the proof, we
mention here that if one equal sign holds the other must also
hold, and that can only happen when all $n$ numbers $a_1 , \ldots
, a_n$ are equal.  We will not prove this here, but the reader may
want to verify this in the cases $n=2$ and $n=3$.  We will proceed
by using mathematical induction to first prove that, $\frac{a_1 +
\ldots + a_n}{n} \geq \sqrt[n]{a_1 \ldots a_n}$, for every natural
number $n$ and all positive reals $a_1, \ldots , a_n$.   Even
though this trivially holds true for $n=1$, we will use as our
starting or base value, $n =2$.  So we first prove that $\frac{a_1
+ a_2}{2} \geq \sqrt{a_1a_2}$ holds true for any two positive
reals.  Since $a_1$ and $a_2$ are both positive, the square roots
$\sqrt{a_1}$ and $\sqrt{a_2}$ are both positive real numbers and
$a_1 = (\sqrt{a_1})^2,\ a_2 = (\sqrt{a_2})^2$.  Clearly,

$$\begin{array}{rl} & (\sqrt{a_1}- \sqrt{a_2})^2 \geq 0 \\
\\
\Ra & (\sqrt{a_1})^2 - 2(\sqrt{a_1})(\sqrt{a_2}) + ( \sqrt{a_2})^2 \geq 0\\
\\
\Ra & a_1 - 2\sqrt{a_1a_2} + a_2 \geq 0 \\
\\
\Ra & a_1 + a_2 \geq 2 \cdot \sqrt{a_1a_2}\\
\\
\Ra & \frac{a_1+a_2}{2} \geq \sqrt{a_1a_2},
\end{array}
$$

\noindent so the statement holds true for $n = 2$.

\vspace{.15in}

\noindent{\bf The Inductive Step:}  Assume the statement to hold
true for some natural number $n = k \geq 2$;  and show that this
assumption implies that the statement must also hold true for $n =
k +1$.  So assume,

$$
\begin{array}{rlll}& \frac{a_1+\ldots + a_k}{k} & \geq & \sqrt[k]{a_1 \ldots a_k} \\
\\
\Ra & a_1 + \ldots + a_k & \geq & k \cdot \sqrt[k]{a_1 \ldots a_k} \end{array}
$$

Now we apply the inequality we proved earlier:

$$b^{k+1} > (k+1) \cdot b \cdot a^k - k\cdot a^{k+1};
$$

\noindent If we take $b = \sqrt[k+1]{a_{k+1}}$, where $a_{k+1}$ is a positive real and $a = \sqrt[k(k+1)]{a_1 \ldots a_k}$ we now have,

$$
\begin{array}{rcl}\left( \sqrt[k+1]{a_{k+1}} \right)^{k+1}& > &(k+1) \cdot \sqrt[k+1]{a_{k+1}} \cdot
\left( \sqrt[k(k+1)]{a_1 \ldots a_k} \right)^k - k \cdot \left( \sqrt[k(k+1)]{a_1 \ldots a_k}\right)^{k+1}\\
\\
& \Ra & a_{k+1} > (k+1)\cdot \sqrt[k+1]{a_{k+1}} \cdot \sqrt[k+1]{a_1 \ldots a_k} - k \cdot \sqrt[k]{a_1 \ldots a_k}\\
\\
& \Ra & a_{k+1} + k \cdot \sqrt[k]{a_1 \ldots a_k} > (k+1)\cdot\sqrt[k+1]{a_1 \ldots a_k \cdot a_{k+1}}
\end{array}
$$

\noindent But from the inductive step we know that $a_1 + \ldots + a_k \geq k\cdot \sqrt[k]{a_1 \ldots a_k}$; hence we have,

$$\begin{array}{rcl}a_{k+1} + (a_1 + \ldots + a_k ) & \geq & a_{k+1} + k\cdot \sqrt[k]{a_1 \ldots a_k} \geq (k+1)\cdot \sqrt[k+1]{a_1
\ldots a_k \cdot a_{k+1}}\\
\\
& \Ra & a_1 + \ldots + a_k + a_{k+1} \geq (k+1) \sqrt[k+1]{a_1 \ldots a_k \cdot a_{k+1}}, \end{array}
$$

\noindent and the induction is complete.

Now that we have established the arithmetic-geometric mean
inequality, we prove the geometric-harmonic inequality.  Indeed,
if $n$ is a natural number and $a_1 , \ldots , a_n$ are positive
reals, then so are the real numbers $\frac{1}{a_1} , \ldots ,
\frac{1}{a_n}$.  By applying the already proven
arithmetic-geometric mean inequality we infer that,

$$
\frac{\frac{1}{a_1} + \ldots + \frac{1}{a_n}}{n} \geq \sqrt[n]{\frac{1}{a_1} \ldots \frac{1}{a_n}}
$$

\noindent Multiplying both sides by the product $\left( \frac{n}{\frac{1}{a_1} + \ldots + \frac{1}{a_n}}\right)
\cdot \sqrt[n]{a_1 \ldots a_n}$, we arrive at the desired result:

$$
\sqrt[n]{a_1 \ldots a_n} \geq \frac{n}{\frac{1}{a_1} +\ldots + \frac{1}{a_n}}.
$$

\noindent This concludes the proof of the theorem. \hfill $\square$

\section{A collection of 21 problems}

\begin{enumerate}
\item[P1.]  Determine the difference of each arithmetic
progression whose first term is $\frac{1}{5}$; and with subsequent
terms (but not necessarily consecutive) the rational numbers
$\frac{1}{4},\ \frac{1}{3},\ \frac{1}{2}$.

\noindent {\bf Solution:}  Let $k,m,n$ be natural numbers with $k
< m < n$ such that $a_k = \frac{1}{4},\ a_m = \frac{1}{3},$ and
$a_n = \frac{1}{2}$.  And, of course, $a_1 = \frac{1}{5}$ is the
first term; $a_1 = \frac{1}{5} , \ldots , a_k = \frac{1}{4},
\ldots , a_m = \frac{1}{3}, \ldots , a_n = \frac{1}{2} , \ldots$ .
By Theorem 1(i) we must have,

\vspace{.15in}

$
\left. \begin{array}{l} \frac{1}{4} = a_k = \frac{1}{5} + (k-1) d\\
\\
\frac{1}{3} = a_m = \frac{1}{5} + (m-1) d\\
\\
\frac{1}{2} = a_n + \frac{1}{5} + (n-1)d \end{array} \right\}$; \begin{tabular}{l} where $d$ is the difference \\
of the arithmetic progression. \end{tabular}

\vspace{.15in}

\noindent Obviously, $d \neq 0$; the three equations yield,

\vspace{.15in}

$\left.\begin{array}{l} (k-1) d = \frac{1}{4} - \frac{1}{5} = \frac{1}{20}\\
\\
(m-1)d = \frac{1}{3} - \frac{1}{5} = \frac{2}{15}\\
\\
(n-1) d = \frac{1}{2} - \frac{1}{5} = \frac{3}{10}\end{array} \right\}$\  \begin{tabular}{ll} (1) & Also, it is clear that $1 < k$;\\
\\
(2) & so that $ 1 < k < m < n$.\\
\\
(3) & \end{tabular}

\vspace{.15in}

\noindent Dividing (1) with (2) member-wise gives

$\frac{k-1}{m-1} = \frac{3}{8},\ \Ra 8(k-1) = 3(m-1)$ (4)

\noindent Dividing (2) with (3) member-wise implies

 $\frac{m-1}{n-1} = \frac{4}{9} \Ra 9(m-1) = 4(n-1)$ (5)

\noindent Dividing (1) with (3) member-wise produces

$\frac{k-1}{n-1} = \frac{1}{6} \Ra 6(k-1) = n-1$  (6)

\vspace{.15in}

According to Equation (4), 3 must be a divisor of $k-1$ and $8$
must be a divisor of $m-1$; if we put $k-1 = 3t;\ k=3t +1$, where
$t$ is a natural number (since $k>1$), then (4) implies $8t = m-1
\Ra m = 8t +1$

Going to equation (5) and substituting for $m-1 = 8t$, we obtain,

$$18t = n-1 \Ra n= 18t + 1.$$

Checking equation (6) we see that $6(3t) = 18t$, which is true for
all nonnegative integer values of $t$. In conclusion we have the
following formulas for $k,\ m,$ and $n$:

$$
k=3t+1,\ m = 8t+1,\ n = 18t + 1;\ t \in {\mathbb N}; \ t = 1,2,\ldots
$$

\noindent We can now calculate $d$ in terms of $t$ from any of the equations (1), (2), or (3):

\noindent From (1), $(k-1)d = \frac{1}{20} \Ra 3t\cdot d =
\frac{1}{20} \Ra$ \fbox{$d = \frac{1}{60t}$.}  We see that this
problem has infinitely many solutions:  there are infinitely many
(infinite) arithmetic progressions that satisfy the conditions of
the problem.  For each positive integer of value of $t$, a new
such arithmetic progression is determined.  For example, for $t =
1$ we have $d=\frac{1}{60},\ k = 4,\ m = 9,\ n = 19$.  We have the
progression,

$$
a_1 = \frac{1}{5},\ldots , a_4 = \frac{1}{4},\ldots \ldots , a_9 = \frac{1}{3}, \ldots \ldots , a_{19} = \frac{1}{2}, \ldots
$$

\item[P2.] Determine the arithmetic progressions (by finding the
first term $a_1$ and difference $d$) whose first term is $a_1 =
5$, whose difference $d$ is an integer, and which contains the
numbers $57$ and $113$ among their terms.

\vspace{.15in}

\noindent {\bf Solution:}  We have $a_1 = 5, \ a_m= 57,\ a_n =
113$ for some natural numbers $m$ and $n$ with $1 < m < n$.  We
have $57 = 5 + (m-1)d$ and $113 = 5 + (n-1)d$;  $(m-1)d = 52$ and
$(n-1)d = 108$;  the last two conditions say that $d$ is a common
divisor of $52$ and $108$;  thus \fbox{$d = 1,2,\ {\rm or}\ 4$}
are the only possible values.  A quick computation shows that for
$d = 1$, we have $m=53$, and $n = 109$; for $d = 2$, we have $m =
27$ and $n = 55$; and for $d = 4, \ m = 14$ and $n = 28$.  In
conclusion there are exactly three arithmetic progressions
satisfying the conditions of this exercise; they have first term
$a_1 = 5$ and their differences $d$ are $d = 1,2,$ and $4$
respectively.

\item[P3.]  Find the sum of all three-digit natural numbers $k$
which are such that the remainder of the divisions of $k$ with
$18$ and of $k$ with $30$,  is equal to $7$.

\noindent{\bf Solution:} Any natural number divisible by both $18$
and $30$, must be divisible by their least common multiple which
is $90$.  Thus if $k$ is any natural number satisfying the
condition of the exercise, then the number $k - 7$ must be
divisible by both $18$ and $90$ and therefore $k-7$ must be
divisible by $90$; so that $k - 7 = 90t$,  for some nonnegative
integer $t$; thus the three-digit numbers of the form $k = 90t+7$
are precisely the numbers we seek to find.  These numbers are
terms in an infinite arithmetic progression whose first term is
$a_1 =7$ and whose difference is $d = 70:\ a_1 = 7,\ a_2 = 7+90,\
a_3 = 7 + 2 \cdot (90), \ldots , a_{t+1}  = 7 + 90t,\ldots$ .

A quick check shows that the first such three-digit number in the
above arithmetic progression is $a_3 = 7 + 90(2) = 187$ (obtained
by setting $t = 2$) and the last such three-digit number in the
above progression is $a_{12} = 7 + 90(11) = 997$ (obtained by
putting $t = 11$ in the formula $a_{t+1} = 7 + 90t$).  Thus, we
seek to find the sum, $a_3 + a_4 + \ldots + a_{11} + a_{12}$.  We
can use either of the two formulas developed in Example 2 (after
example 1 which in turn is located below the proof of Theorem 2).

Since we know the first and last terms of the sum at hand, namely $a_3$, it is easier to use the first formula in Example 2:

$$
\begin{array}{rcl}a_m + a_{m+1} + \ldots + a_{n-1} + a_n & = & \frac{(n-m+1)(a_m+a_n)}{2}
\end{array}$$

\noindent In our case $m = 3,\ n = 12, \ a_m = a_3 = 187$, and $a_n = a_{12} = 997$.  Thus

$$\begin{array}{rcl} a_3 +a_4 + \ldots + a_{11} + a_{12} & = & \frac{(12-3+1)\cdot(187+997)}{2}\\
\\
&  = & \frac{10}{2} \cdot (1184) = 5 \cdot (1184) = 5920.\end{array}
$$

\item[P4.]  Let $a_1, a_2, \ldots , a_n,\ldots$, be an arithmetic
progression with first term $a_1$ and positive difference $d$; and
$M$ a natural number, such that $a_1 \leq M$.  Show that the
number of terms of the arithmetic progression that do not exceed
$M$, is equal to $\left[\!\left[ \frac{M-a_1}{d}\right]\!\right] +
1$, where $\left[\!\left[ \frac{M-a_1}{d} \right]\!\right]$ stands
for the integer part of the real number $\frac{M-a_1}{d}$.

\vspace{.15in}

\noindent{\bf Solution:}  If, among the terms of the arithmetic
progression, $a_n$ is the largest term which does not exceed $M$,
then $a_n \leq M$ and $a_{\ell} > M$, for all natural number
$\ell$ greater than $n$;  $\ell = n+1,n+2,\ldots$ .  But $a_n =
a_1 + (n-1)d$; so that $a_1 + (n-1)d\leq M\Ra (n-1)d \leq M - a_1
\Ra n-1 \leq \frac{M-a_1}{d}$ since $d > 0$.  Since, by
definition, $\left[\!\left[ \frac{M-a_1}{d} \right]\!\right]$  is
the greatest integer not exceeding $\frac{M-a_1}{d}$ and since
$n-1$ does not exceed $\frac{M-a_1}{d}$, we conclude that $n - 1
\leq \left[\!\left[ \frac{M-a_1}{d} \right]\!\right] \Ra n \leq
\left[\!\left[ \frac{M-a_1}{d} \right]\!\right]  +1$. But $n$ is a
natural number, that is, a positive integer, and so must be the
integer $N = \left[\!\left[ \frac{M-a_1}{d} \right]\!\right] +1$
Since $a_n$ was assumed to be the largest term such that $a_n \leq
M$, it follows that $n$ must equal $N$;  because the term $a_N$ is
actually the largest term not exceeding $M$ (note that if $n < N$,
then $a_n < a_N$, since the progression is increasing in view of
the fact that $d > 0$).  Indeed, if $N = \left[\!\left[
\frac{M-a_1}{d}\right]\!\right] + 1$, then by the definition of
the integer part of a real number we must have $N - 1 \leq
\frac{M-a_1}{d}<N$.  Multiplying by $d > 0$ yields $d(N-1) \leq M
- a_1 \Ra a_1 + d(N-1) \leq M \Ra a_N \leq M$.

In conclusion we see that the terms $a_1, \ldots , a_N$ are
precisely the terms not exceeding $\left[\!\left[
\frac{M-a_1}{d}\right]\!\right] +1$; therefore there are exactly
$\left[\!\left[ \frac{M-a_1}{d} \right]\!\right] + 1$ terms not
exceeding $M$.

\item[P5.]  Apply the previous problem P4 to find the value of the
sum of all natural numbers $k$ not exceeding $1,000$, and which
are such that the remainder of the division of $k^2$ with $17$ is
equal to $9$.

\vspace{.15in}

\noindent{\bf Solution:}  First, we divide those numbers $k$ into
two disjoint classes or groups.  If $q$  is the quotient of the
division of $k^2$ with $17$, and with remainder $9$, we must have,

$$
k^2 = 17q + 9 \Lra (k-3)(k+3) = 17q,
$$

\noindent but $17$ is a prime number and as such it must divide at
least one of the two factors $k-3$ and $k+3$; but it cannot divide
both.  Why?  Because for any value of the natural number $k$, it
is easy to see that the greatest common divisor of $k-3$ and $k+3$
is either equal to $1,2$, or $6$. Thus, we must have either $k - 3
= 17n$ or $k+3 = 17m$; either $k = 17n+3$ or

$$\begin{array}{rcl} k=17 m-3 & = &  17(m-1) + 14\\
& = & 17 \cdot \ell + 14 \end{array}
$$

\noindent (here we have set $m-1 = \ell$). The number $n$ is a
nonnegative integer and the number $\ell$ is also a nonnegative
integer.  So the two disjoint classes of the natural numbers $k$
are,

$$
\begin{array}{rrcl}& k & = & 3, 20, 37, 54, \ldots\\
\\
{\rm and} & k & = & 14, 31, 48, 65, \ldots \end{array}
$$

\noindent Next, we find how many numbers $k$ in each class do not
exceed $M = 10,000$.  Here, we are dealing with two arithmetic
progressions:  the first being $3,20,37,54,\ldots ,$ having first
term $a_1 = 3$ and difference $d = 17$. The second arithmetic
progression has first term $b_1 = 14 $ and the same difference $d
= 17$.

According to the previous practice problem, P4, there are exactly
$N_1 = \left[\!\left[ \frac{M-a_1}{d} \right]\!\right] + 1 =
\left[\!\left[ \frac{1000 - 3}{17}\right]\!\right] + 1 =
\left[\!\left[ \frac{997}{17} \right]\!\right] + 1 = 58 + 1 = 59$
terms of the first arithmetic progression not exceeding $1000$
(also, recall from Chapter 6 that $\left[\!\left[
\frac{997}{17}\right]\!\right]$ is really none other than the
quotient of the division of $997$ with $17$).

Again, applying problem P4 to the second arithmetic progression,
we see that there are $N_2 = \left[\!\left[ \frac{M-b_1}{d}
\right]\!\right] + 1 = \left[\!\left[ \frac{1000-14}{17}
\right]\!\right] + 1 = \left[\!\left[ \frac{986}{17}
\right]\!\right] + 1 = 58+1 = 59$.

\noindent Finally, we must find the two sums:

$$
\begin{array}{rcl} S_{N_1} & = & a_1 + \ldots + a_{N_1} = \frac{N_1 \cdot (a_1 + a_{N_1})}{2} = \frac{N_1 \cdot \left[ 2a_1 + (N_1 -1)d\right]}{2} \\
\\
& = & \frac{59\cdot \left[2(3) + (59-1) \cdot 17 \right]}{2} = \frac{59 \cdot \left[ 6 + (58)(17)\right]}{2}
\end{array}
$$

\noindent and

$$
\begin{array}{rcl}
S_{N_2} & = & b_1 \ldots + b_{N_2} = \frac{N_2 \cdot \left[ 2b_1 + (N_2 - 1)d\right]}{2} \\
\\
& = & \frac{59 \cdot \left[2(14)+ (59-1)17\right]}{2} = \frac{59 \cdot \left[ 28 + (58)(17)\right]}{2}
\end{array}
$$

\noindent Hence,

$$
\begin{array}{rcl}
S_{N_1} + S_{N_2} & = & \frac{59\cdot \left[ 6 + 28 + 2(58)(17)\right]}{2}\\
\\
& = & \frac{59\left[34 + 1972\right]}{2} = \frac{59 \cdot (2006)}{2} = 59 \cdot (1003) = 59,177.
\end{array}
$$

\item[P6.]  If $S_n,\ S_{2n},\ S_{3n}$, are the sums of the first
$n,\ 2n,\ 3n$ terms of an arithmetic progression, find the
relation or equation between the three sums.

\vspace{.15in}

\noindent{\bf Solution:}  We have $S_n = \frac{n\cdot \left[ a_1 +
(n-1)d\right]}{2}$, $S_{2n}= \frac{2n\cdot \left[ a_1 +
(2n-1)d\right]}{2}$, \linebreak and $S_{3n} = \frac{3n\cdot
\left[a_1 + (3n-1)d\right]}{2}$.

We can write
$$
\begin{array}{rcl}
S_{2n} & = & \frac{2n\cdot \left[ 2a_1 + 2(n-1)d+(d-a_1)\right]}{2}\ {\rm and}\\
\\
S_{3n} & = & \frac{3n \cdot \left[3a_1 + 3(n-1)d + (2d-2a_1)\right]}{2}.
\end{array}
$$
\noindent So that,

\hspace*{1.0in}$S_{2n}  =  \frac{2n\cdot 2 \cdot \left[ a_1 + (n-1)d\right]}{2} + \frac{2n\cdot (d-a_1)}{2}$\hfill (1)

\vspace{.15in}

\noindent and

\hspace*{1.0in}$S_{3n} = \frac{3n\cdot 3 \cdot \left[ a_1 + (n-1)d\right]}{2} + \frac{3n \cdot 2\cdot (d - a_1)}{2}$ \hfill (2)

\vspace{.15in}

To eliminate the product $n \cdot (d-a_1)$ in equations (1) and (2) just consider $3S_{2n}-S_{3n}$: equations (1) and (2) imply,

$$
\begin{array}{rcl}
3S_{2n} - S_{3n} & = & \frac{3 \cdot 2n \cdot 2 \cdot \left[ a_1 + (n-1)d\right]}{2} - \frac{3n\cdot 3 \cdot \left[ a_1 + (n-1)d\right]}{2} \\
\\
& & + \underset{0}{\underbrace{\frac{3 \cdot 2n \cdot (d-a_1)}{2} - \frac{3n \cdot 2 \cdot (d - a_1)}{2}}} \\
\\
\Ra 3S_{2n} - S_{3n} & = & \frac{3n\cdot \left[ a_1 + (n-1)d\right]}{2}
\end{array}
$$

\noindent but $S_n = \frac{n\cdot \left[ a_1 + (n-1)d\right]}{2}$; hence the last equation yields

$$
\begin{array}{rl}
& 3S_{2n} - S_{3n}  = 3 \cdot S_n\\
\\
\Ra & \fbox{$3S_{2n} = 3S_n + S_{3n}$};\\
\\
{\rm or} & 3(S_{2n} - S_n) = S_{3n}
\end{array}
$$

\item[P7.]  If the first term of an arithmetic progression is
equal to some real number $a$, and the sum of the first $m$ terms
is equal to zero, show that the sum of the next $n$ terms must
equal to $\frac{a\cdot m(m+n)}{1-m}$; here, we assume that $m$ and
$n$ are natural numbers with $m > 1$

\vspace{.15in}

\noindent{\bf Solution:}  We have $a_1 + \ldots + a_m = 0 =
\frac{m\cdot \left[2a_1 + d(m-1)\right]}{2} \Ra$ (since $m > 1$)
$2a_1 + d(m-1) = 0 \Ra d = \frac{-2a_1}{m-1} = \frac{2a_1}{1-m} =
\frac{2a}{1-m}$.  Consider the sum of the next $n$ terms

$$
\begin{array}{rcl}
a_{m+1}+\ldots + a_{m+n} & = & \frac{n\cdot (a_{m+1} + a_{m+n})}{2} ;\\
\\
a_{m+1} +\ldots + a_{m+n} & = & \frac{n \cdot \left[(a_1 + md) + (a_1 + (m+n-1)d)\right]}{2};\\
\\
a_{m+1}+ \ldots + a_{m+n} & = & \frac{n \cdot \left[ 2a_1 + (2m+n-1)d\right]}{2}
\end{array}
$$

Now substitute for $d = \frac{2a}{1-m}$: (and of course, $a = a_1$)

$$
\begin{array}{rcl}
a_{m+1} + \ldots + a_{m+n} & = & \frac{n[2a +(2m+n-1)\cdot \frac{2a}{1-m}]}{2};\\
\\
a_{m+1} + \ldots + a_{m+n} & = & \frac{n\cdot 2a[(1-m)+(2m+n-1)]}{2(1-m)} ;\\
\\
a_{m+1} + \ldots + a_{m+n} & = & \frac{2an[1-m+2m+n-1]}{2(1-m)} = \frac{a\cdot n\cdot (m+n)}{1-m}
\end{array}
$$

\item[P8.]  Suppose that the sum of the $m$ first terms of an
arithmetic progression is $n$; and thesum of the first $n$ terms
is equal to $m$.  Furthermore, suppose that the first term is
$\alpha$ and the difference is $\beta$, where $\alpha$ and $\beta$
are given real numbers.  Also, assume $m \neq n$ and $\beta \neq
0$.

\begin{enumerate}
\item[(a)] Find the sum of the first $(m + n)$ in terms of the
constants $\alpha$ and $\beta$ only.

\item[(b)] Express the integer $mn$ and the difference $(m-n)$ in
terms of $\alpha$ and $\beta$.

\item[(c)]  Drop the assumption that $m \neq n$, and suppose that
both $\alpha$ and $\beta$ are integers.  Describe all such
arithmetic progressions.
\end{enumerate}

\noindent{\bf Solution:}
\begin{enumerate}
\item[(a)] We have $a_1 + \ldots + a_m = n$ and $a_1 + \ldots + a_n = m$;

$$
\frac{m \cdot[2\alpha +(m-1)\beta]}{2} = n\ \ {\rm and}\ \ \frac{n\cdot[2\alpha + (n-1)\beta]}{2} = m,
$$

\noindent since $a_1 = \alpha$ and $d = \beta$.

Subtracting the second equation from the first one to obtain,

$$
\begin{array}{rcl}
2\alpha \cdot (m-n) & + & \beta \cdot [m(m-1)-n(n-1)] = 2n -2m;\\
2\alpha \cdot (m-n) & + & \beta \cdot [(m^2 - n^2) - (m-n)] + 2(m-n) = 0;\\
2\alpha \cdot (m-n) & + & \beta \cdot [(m-n)(m+n) - (m-n)] + 2(m-n) = 0;\\
2\alpha \cdot (m-n) & + & \beta \cdot (m-n) \cdot [m+n-1]+2(m-n)=0;
\end{array}
$$

\noindent $(m-n) \cdot [2\alpha + \beta(m+n-1) + 2] = 0$; but $m-n\neq 0$, since $m \neq n$ by the hypothesis of the problem.  Thus,

$$
\begin{array}{ll}
& 2 \alpha + \beta\cdot (m+n-1) + 2 = 0 \Ra \beta (m+n-1) = -2(1 + a)\\
\\
\Ra & m + n - 1 = \frac{-2(1+a)}{\beta} \Ra m + n = 1 -\frac{2(1+\alpha)}{\beta} = \frac{\beta - 2\alpha - 2}{\beta}.
\end{array}
$$

\noindent Now, we compute the sum $a_1 + \ldots + a_{m+n} = \frac{(m+n)\cdot[2\alpha +(m+n-1)\beta]}{2}$

$$
\begin{array}{rl}
\Ra & a_1+\ldots +a_{m+n} = \frac{\left(\frac{\beta - 2\alpha -2}{\beta}\right)\cdot \left[ 2 \alpha
\left(\frac{\beta-2\alpha - 2}{\beta} \right) \cdot \beta \right]}{2};\\
\\
& a_1 + \ldots + a_{m+n} = \fbox{$\frac{\left(\beta - 2\alpha -2\right) \cdot \left(\beta - 2\right)}{2\beta}$}
\end{array}
$$

\item[(b)]  If we multiply the equations $\frac{m\cdot[2\alpha
+(m-1)\beta]}{2} = n$ and $\frac{n\cdot[2 \alpha + (n-1)\beta]}{2}
= m$ member-wise we obtain, $\frac{m\cdot n\cdot
[2\alpha+(n-1)\beta][2\alpha +(m-1)\beta]}{4} = mn$ and since $mn
\neq 0$, we arrive at

$$
\begin{array}{rl} & [2\alpha + (n-1)\beta]\cdot[2\alpha +(m-1)\beta] = 4\\
\\
\Ra & 4\alpha^2 + 2\alpha \beta \cdot (m-1+n-1) +(n-1)(m-1)\beta^2 = 4\\
\\
\Ra & 4\alpha^2 + 2\alpha \beta \cdot (m+n) - 4\alpha \beta +nm\beta^2 -(n+m)\beta^2 + \beta^2 = 4;\\
\\
& (2\alpha - \beta )^2 + (m+n)\cdot (2\alpha \beta - \beta^2) + nm\beta^2 = 4. \end{array}
$$

\noindent Now let us substitute for $m+n = \frac{\beta - 2\alpha -2}{\beta}$ (from part (a)) in the last equation above; we have,

$$
\begin{array}{rl}
& (2\alpha-\beta)^2 + \left( \frac{\beta - 2\alpha -2}{\beta}\right) \cdot \beta \cdot (2\alpha -\beta) + nm\beta^2 = 4\\
\\
\Ra & (2\alpha - \beta)^2 + (\beta - 2\alpha-2)(2\alpha - \beta) + nm\beta^2 = 4\\
\\
\Ra & 4\alpha ^2 - 4\alpha \beta + \beta^2 + 2\alpha \beta - \beta^2 - 4 \alpha^2 + 4 \alpha \beta -4\alpha + 2\beta +nm\beta^2 = 4\\
\\
\Ra & nm\beta^2 + 2\alpha \beta - 4\alpha + 2\beta = 4 \Ra nm\beta^2 = 4 - 2\alpha \beta + 4\alpha -2\beta\\
\\
\Ra & \fbox{$nm = \frac{2\cdot(2-\alpha \beta + 2\alpha - \beta)}{\beta^2}$}
\end{array}
$$

Finally, from the identity $(m-n)^2 = (m+n)^2 - 4nm$, it follows that

$$
\begin{array}{rl}
&(m-n)^2 = \left(\frac{\beta - 2\alpha -2}{\beta}\right)^2 - \frac{8(2-\alpha\beta + 2\alpha - \beta)}{\beta^2}\\
\\
\Ra & (m-n)^2 = \frac{\beta^2+4\alpha^2 + 4 - 4 \alpha \beta - 4\beta + 8 \alpha-16 + 8\alpha \beta - 16\alpha + 8\beta}{\beta^2}\\
\\
& (m - n)^2 = \frac{\beta^2 + 4\alpha^2 -12 + 4 \alpha \beta +4\beta -8\alpha}{\beta^2};\\
\\
& |m - n| = \frac{\sqrt{\beta^2 + 4\alpha^2 - 12 + 4 \alpha \beta + 4\beta -8\alpha}}{|\beta |}\\
&  = \frac{\sqrt{(2\alpha + \beta)^2 -12+4\beta - 8 \alpha}}{|\beta|};\\
\\
& \fbox{$m-n = \pm \frac{\sqrt{(2\alpha + \beta)^2 - 12 + 4\beta - 8\alpha}}{|\beta |}$}
\end{array}
$$

\noindent the choice of the sign depending on whether $m > n$ or $m < n$ respectively.  Also note, that a necessary condition that must hold here is

$$
(2\alpha + \beta)^2 -12 + 4\beta -8\alpha > 0.
$$

\item[(c)] Now consider $\dfrac{m[2\alpha +(m-1)\beta]}{2} = n$
and $\dfrac{n[2\alpha +(n-1)\beta]}{2} = m$, with $\alpha$ and
$\beta$ being integers.  There are four cases.

\vspace{.15in}

\noindent {\bf Case 1:}  Suppose that $m$ and $n$ are odd.  Then
we see that $m\mid n$ and $n \mid m$, which implies $m = n$ (since
$m,n$ are positive integers; if they are divisors of each other,
they must be equal).  We obtain,

$$
2\alpha + (n-1) \beta = 2 \Leftrightarrow n =
\dfrac{\beta+2-2\alpha}{\beta} = 1 + \dfrac{2(1-\alpha)}{\beta};\
\beta \mid 2(1-\alpha).
$$

If $\beta$ is odd, it must be a divisor of $1 - \alpha$.  Put
$1-\alpha = \beta \rho$ and so $n = 1+2\rho$, with $\rho$ being a
positive integer.  So, the solution is

\vspace{.15in}

\begin{center}\fbox{\parbox{3.5in}{$m = n = 1 + 2\rho,\ \ \alpha = 1 -
\beta \rho,\ \ \rho \in \Z^+, \ \ \beta \in \Z$}}\end{center}

\vspace{.15in}

If $\beta $ is even, set $\beta = 2B$.  We obtain $1 - \alpha -
B\rho$, for some odd integer $\rho \geq 1$.  The solution is

\vspace{.15in}

\fbox{\parbox{4.5in}{$m = n = 1 + \rho,\ \ \alpha = 1- B\rho, \ \
\beta = 2B,\ \ \rho\ {\rm an \ odd\ positive\ integer}$.}}

\vspace{.15in}

\noindent {\bf Case 2:}  Suppose that $m $ is even, $n$ is odd;
put $m = 2k$.  We obtain

$$
k \left[ 2 \alpha + (2k-1)\beta\right] = n\ {\rm and}\ n\left[ 2
\alpha + (n-1) \beta \right] = 4k.
$$

Since $n$ is odd, $n $ must be a divisor of $k$ and since $k $ is
also a divisor of $n$, we conclude that since $n$ and $k$ are
positive, we must have $n = k$.  So, $2 \alpha + (2n-1) \beta =1$
and $2 \alpha + (n-1) \beta = 4$. From which we obtain $n \beta =
-3 \Leftrightarrow (n=1\ {\rm and }\ \beta = -3)$ or $(n=3\ {\rm
and}\ \beta = 01$).

The solution is

\vspace{.15in}

\begin{center}\fbox{\parbox{2.75in}{$\begin{array}{rl}& n = 1,\ \beta = -3,\
m =
2,\ \alpha = 2 \\
{\rm or} & n=3,\ \beta = -1,\ m = 6,\ \alpha =
3\end{array}$}}\end{center}

\vspace{.15in}

\noindent {\bf Case 3:}  $m$ odd and $n$ even.  This is exactly
analogous to the previous case.  One obtains the solutions (just
switch $m$ and $n$)

\vspace{.15in}

\begin{center}\fbox{\parbox{2.5in}{$\begin{array}{l} m=1,\ \beta = -3,\ n =
2,\ \alpha = 2\\ m = 3,\ \beta = -1,\ n = 6,\ \alpha = 3
\end{array}$}}\end{center}

\vspace{.15in}

\noindent {\bf Case 4:}  Assume $m$ and $n$ to be both even.  Set
$m = 2^e_{m_1}, n = 2^f_{n_1}$, where $e,f$ are positive integers
and $m_1,n_1$ are odd positive integers. Since $n-1$ and $m-1$ are
odd, by inspection we see that $\beta$ must be even.  We have,

$$
\left\{\begin{array}{rl}& 2^e \cdot m_1 \cdot \left[ 2\alpha
+\left(2^3_{m_1} -1\right)\cdot \beta\right] = 2^{f+1} \cdot n_1 \\
\\
{\rm and }& 2^f \cdot n_1 \cdot \left[ 2\alpha + \left[ 2\alpha +
\left(2^f_{n_1}-1\right) \cdot \beta \right]\right] = 2^{e+1}
\cdot m_1.
\end{array}\right.$$

We see that the left-hand side of the first equation is divisible
by a power of 2 which is at least $2^{e+1}$; and the left-hand
side of the equation is divisible by at least $2^{f+1}$.

This then implies that $e+1 \leq f+1$ and $f+1 \leq e+1$.  Hence
$e = f$.  Consequently,

$$\begin{array}{rcll}
m_1 \left[ 2 \alpha + \left(2^e_{m_1}-1\right) \beta \right]
&=& 2_{n_1} & {\rm and}\\
\ n_1 \left[ 2 \alpha + \left(2^e_{n_1} -1\right) \beta \right]&
=& 2m_1 &
\end{array}
$$

Let $\beta = 2k$.  By cancelling the factor 2 from both sides of
the two equations, we infer that $m_1$ is a divisor of $n_1$ and
$n_1$ a divisor of $m_1$.  Thus $m_1 = n_1$.

The solution is

\vspace{.15in}

\begin{center}\fbox{\parbox{1.75in}{$\begin{array}{l} \alpha = 1 - \left(
2^e
\cdot n_1-1 \right) k \\
\\
\beta = 2k \\
\\
m = 2^e_{n_1} = n\end{array}$ }},\end{center}

\noindent where $k$ is an arbitrary integer, $e$ is a positive
integer, and $n_1$ can be any odd positive integer.
\end{enumerate}

\item[P9.] Prove that if the real numbers $\alpha,\beta, \gamma, \delta$ are successive terms of a harmonic progression, then

$$
3(\beta-\alpha) (\delta - \gamma) = (\gamma - \beta)(\delta - \alpha).
$$

\vspace{.15in}

\noindent{\bf Solution:}  Since $\alpha, \beta, \gamma,\delta$ are members of a harmonic progression they must all be nonzero;
$\alpha \beta \gamma \delta \neq 0$.  Thus

$$
3(\beta - \alpha)(\delta - \gamma) = ( \gamma - \beta)(\delta - \alpha)
$$

\noindent is equivalent to

$$
\frac{3(\beta-\alpha)(\delta - \gamma)}{\alpha \beta \gamma \delta} = \frac{(\gamma - \beta)(\delta - \alpha)}{\alpha\beta\gamma\delta}
$$

$$
\begin{array}{rl}
\Lra & 3\cdot \left( \frac{\beta - \alpha}{\beta \alpha}\right) \cdot \left( \frac{\delta - \gamma}{\delta \gamma}\right) =
\left( \frac{\gamma - \beta}{\gamma \beta}\right) \cdot \left( \frac{\delta - \alpha}{\alpha \delta}\right)\\
\\
\Lra & 3 \cdot \left( \frac{1}{\alpha} - \frac{1}{\beta} \right) \cdot \left(\frac{1}{\gamma} - \frac{1}{\delta}\right) =
\left(\frac{1}{\beta}-\frac{1}{\gamma}\right) \cdot \left( \frac{1}{\alpha} - \frac{1}{\delta} \right)
\end{array}
$$

\noindent By definition, since $\alpha,\beta,\gamma,\delta$ are
consecutive terms of a harmonic progression; the numbers
$\frac{1}{\alpha},\frac{1}{\beta},\frac{1}{\gamma},\frac{1}{\delta}$
must be successive terms of an arithmetic progression with
difference $d$; and $\frac{1}{\alpha} - \frac{1}{\beta} = -d,\
\frac{1}{\gamma}-\frac{1}{\delta} = -d$, \linebreak
$\frac{1}{\beta} - \frac{1}{\gamma} = -d$, and $ \frac{1}{\alpha}
- \frac{1}{\delta} = -3d$ (since $\frac{1}{\delta} =
\frac{1}{\gamma} + d = \frac{1}{\beta} + 2d = \frac{1}{\alpha} +
3d$).  Thus the above statement we want to prove is equivalent to

$$3 \cdot (-3)\cdot(-d) = (-d)\cdot(-3d) \Lra 3d^2 = 3d^2$$

\noindent which is true.

\item[P10.]  Suppose that $m$ and $n$ are fixed natural numbers
such that the $m$th term $a_m$ in a harmonic progression is equal
to $n$;  and the $n$th term $a_n$ is equal to $m$.  We assume $m
\neq n$.

\begin{enumerate}
\item[(a)] Find the $(m+n)$th term $a_{m+n}$ in terms of $m$ and $n$ .
\item[(b)] Determine the general $k$th term $a_k$ in terms of $k,m$, and $n$.
\end{enumerate}

\noindent {\bf Solution:}
\begin{enumerate}
\item[(a)] Both $\frac{1}{a_m}$ and  are the $m$th and $n$th terms
respectively of an arithmetic progression with first term
$\frac{1}{a_1}$ and difference $d$; so that $\frac{1}{a_m} =
\frac{1}{a_1} + (m-1)d$ and $\frac{1}{a_n} = \frac{1}{a_1} +
(n-1)d$. Subtracting the second equation from the first and using
the fact that $a_m = n$ and $a_n=m$ we obtain, $\frac{1}{n} -
\frac{1}{m} = (m-n)d \Ra \frac{m-n}{nm} = (m-n)d$; but $m-n \neq
0$; cancelling the factor $(m-n)$ from both sides, gives
\fbox{$\frac{1}{mn} = d$}.  Thus from the first equation ,
$\frac{1}{n} = \frac{1}{a_1} + (m-1)\cdot \frac{1}{mn}\Ra
\frac{1}{n} -\frac{(m-1)}{mn}= \frac{1}{a_1} \Ra
\frac{m-(m-1)}{mn} = \frac{1}{a_1} ;\ \frac{1}{mn} = \frac{1}{a_1}
\Ra  \fbox{$a_1 = mn$}$.  Therefore, $\frac{1}{a_{m+n}} =
\frac{1}{a_1} + (m+n-1)d \Ra \frac{1}{a_{m+n}} = \frac{1}{mn} +
\frac{m+n-1}{mn} \Ra \fbox{${a}_{m+n} = \frac{mn}{m+n}$}$.

\item[(b)] We have $\frac{1}{a_k} = \frac{1}{a_1} + (k-1)d \Ra
\frac{1}{a_k} = \frac{1}{mn}+ \frac{(k-1)}{mn} = \frac{k}{mn} \Ra
\fbox{$a_k = \frac{mn}{k}$}$ .
\end{enumerate}

\item[P11.]  Use mathematical induction to prove that if $a_1,a_2,
\ldots, a_n$, with $n \geq 3$, are the first $n$ terms of a
harmonic progression, then $(n-1)a_1a_n = a_1a_2+a_2a_3 + \ldots +
a_{n-1}a_n$.

\vspace{.15in}

\noindent{\bf Solution:}  For $n=3$ the statement is $2a_1 a_3 =
a_1a_2+a_2a_3 \Lra a_2 \cdot (a_1+a_3)=2a_1a_3$; but $a_1,a_2,a_3$
are all nonzero since they are the first three terms of a harmonic
progression.  Thus, the last equation is equivalent to
$\frac{2}{a_2} = \frac{a_1+a_d}{a_1a_3} \Lra \frac{2}{a_2} =
\frac{1}{a_3} + \frac{1}{a_1}$ which is true, because
$\frac{1}{a_1},\frac{1}{a_2},\frac{1}{a_3}$ are the first three
terms of a harmonic expression.

The inductive step:  prove that whenever the statement holds true
for some natural number $n = k \geq 3$, then it must also hold
true for $n = k+1$.  So we assume $(k-1)a_1a_k = a_1a_2+a_2a_3 +
\ldots + a_{k-1}a_k$.  Add $a_ka_{k+1}$ to both sides to obtain,

\vspace{.15in}

$(k-1)a_1a_k + a_ka_{k+1} = a_1a_2+a_2a_3 + \ldots + a_{k-1}a_k + a_ka_{k+1}$ \hfill (1)

\vspace{.15in}

\noindent If we can show that the left-hand side of (1) is equal
to $ka_1a_{k+1}$, the induction process will be complete.  So we
need to show that

\vspace{.15in}

\hspace*{1.0in}$(k-1)a_1a_k+a_ka_{k+1} = k\cdot a_1 \cdot a_{k+1}$ \hfill (2)

\vspace{.15in}

\noindent (dividing both sides of the equation by $a_1\cdot a_k\cdot a_{k+1} \neq 0$)

\vspace{.15in}

\hspace{.15in} $\Lra \frac{(k-1)}{a_{k+1}} + \frac{1}{a_1} = \frac{k}{a_k}.$\hfill (3)

\vspace{.15in}

To prove (3), we can use the fact that $\frac{1}{a_{k+1}}$ and
$\frac{1}{a_k}$ are the $(k+1)$th and $k$th terms of an arithmetic
progression with first term $\frac{1}{a_1}$ and ratio $d$:
$\frac{1}{a_{k+1}} = \frac{1}{a_1} + k \cdot d$ and $\frac{1}{a_k}
= \frac{1}{a_1} + (k-1)d$;  so that, $\frac{k-1}{a_{k+1}} =
\frac{k-1}{a_1} + (k-1)kd$ and $\frac{k}{a_k} = \frac{k}{a_1} +
k(k-1)d$.  Subtracting the second equation from the first yields,
$$
\frac{k-1}{a_{k+1}} - \frac{k}{a_k} = \frac{(k-1)-k}{a_1} \Ra \frac{k-1}{a_{k+1}} + \frac{1}{a_1} = \frac{k}{a_k}
$$
\noindent which establishes (3) and thus equation (2).  The induction is complete since we have show (by combining (1) and (3)).

$$
k\cdot a_1a_{k+1} = a_1a_2 +a_2a_3 +\ldots +a_{k-1}a_k+a_k a_{k+1},
$$

\noindent the statement also holds for $n=k+1$.

\item[P12.]  Find the necessary and sufficient condition that
three natural numbers $m,n$, and $k$ must satisfy, in order that
the positive real numbers $\sqrt{m},\sqrt{n},\sqrt{k}$ be
consecutive terms of a geometric progression.

\vspace{.15in}

\noindent{\bf Solution:}  According to Theorem 7, the three
positive reals will be consecutive terms of an arithmetic
progression if, and only if, $(\sqrt{n})^2 = \sqrt{m}\sqrt{k} \Lra
n = \sqrt{mk} \Lra$ (since both $n$ and $mk$ are positive) $n^2 =
mk$.  Thus, the necessary and sufficient condition is that the
product of $m$ and $k$ be equal to the square of $n$.

\item[P13.]  Show that if $\alpha, \beta, \gamma$ are successive
terms of an arithmetic progression, $\beta,\gamma, \delta$ are
consecutive terms of a geometric progression, and $\gamma, \delta,
\epsilon$ are the successive terms of a harmonic progression, then
either the numbers $\alpha, \gamma, \epsilon$ or the numbers
$\epsilon, \gamma, \alpha$ must be the consecutive terms of a
geometric progression.

\vspace{.15in}

\noindent{\bf Solution:} Since
$\frac{1}{\gamma},\frac{1}{\delta},\frac{1}{\epsilon}$ are by
definition successive terms of an arithmetic progression and the
same holds true for $\alpha , \beta, \gamma$, Theorem 3 tells us
that we must have $2\beta = \alpha + \gamma$ (1) and
$\frac{2}{\delta} = \frac{1}{\gamma}+\frac{1}{\epsilon}$ (2). And
by Theorem 7, we must also have $\gamma^2 = \beta \delta$ (3).
(Note that $\gamma, \delta$, and $\epsilon$ must be nonzero and
thus so must be $\beta$.)

Equation (2) implies $\delta = \frac{2\gamma \epsilon}{\gamma +
\epsilon}$ and equation (1) implies $\beta = \frac{\alpha +
\gamma}{2}$.  Substituting for $\beta$ and $\delta$ in equation
(3) we now have

$$\begin{array}{rl}
& \gamma^2 = \left( \frac{\alpha + \gamma}{2}\right) \cdot \left( \frac{2\gamma \epsilon}{\gamma+\epsilon}\right)\\
\\
 \Ra&  \gamma^2 \cdot (\gamma + \epsilon) = (\alpha + \gamma)\cdot \gamma \epsilon \Ra \gamma^3 + \gamma^2
 \epsilon = \alpha \gamma \epsilon + \gamma^2\epsilon \\
\\
\Ra & \gamma^3 - \alpha \gamma \epsilon = 0 \Ra \gamma (\gamma^2 - \alpha \epsilon) =0\end{array}
$$

\noindent and since $\gamma \neq 0$ we conclude $\gamma^2 - \alpha
\epsilon = 0 \Ra \gamma^2 = \alpha\epsilon$, which, in accordance
with Theorem 7, proves that either $\alpha, \gamma, \epsilon$; or
$\epsilon, \gamma, \alpha$ are consecutive terms in a geometric
progression.

\item[P14.]  Prove that if $\alpha$ is the arithmetic mean of the
numbers $\beta$ and $\gamma$; and $\beta$, nonzero, the geometric
mean of $\alpha$ and $\gamma$, then $\gamma$ must be the harmonic
mean of $\alpha$ and $\beta$.  (Note:  the assumption $\beta \neq
0$, together with the fact that $\beta$ is the geometric  mean of
$\alpha$ and $\gamma$, does imply that both $\alpha$ and $\gamma$
must be nonzero as well.)

\vspace{.15in}

\noindent{\bf Solution:}  From the problems assumptions we must
have $2\alpha = \beta + \gamma$ and $\beta^2 = \alpha \gamma$;
$\beta^2 = \alpha \gamma \Ra 2\beta^2 = 2 \alpha \gamma$;
substituting for $2\alpha = \beta + \gamma$ in the last equation
produces

$$
\begin{array}{rl}
& 2\beta^2 = (\beta+\gamma)\gamma \Ra 2\beta^2 = \beta\gamma + \gamma^2\\
\\
 \Ra  & 2\beta^2 - \gamma^2 - \beta \gamma = 0 \Ra (\beta^2 - \gamma^2) +(\beta^2 - \beta\gamma) = 0\\
\\
 \Ra &  (\beta-\gamma)(\beta+\gamma) + \beta \cdot (\beta - \gamma) = 0 \Ra (\beta - \gamma) \cdot (2\beta + \gamma) = 0.
\end{array}
$$

\noindent If $\beta - \gamma \neq 0$, then the last equation
implies $2\beta + \gamma = 0 \Ra \gamma = -2\beta$; and thus from
$2a=\beta + \gamma$ we obtain $2\alpha = \beta -2\beta$; $2\alpha
= -\beta$;  $\alpha = -\beta/2$.  Now compute, $\frac{2}{\gamma} =
\frac{2}{-2\beta} = - \frac{1}{\beta}$, since $\beta \neq 0$; and
$\frac{1}{\alpha} + \frac{1}{\beta} = \frac{1}{-\frac{\beta}{2}} +
\frac{1}{\beta} = - \frac{2}{\beta} + \frac{1}{\beta} = -
\frac{1}{\beta}$.  Therefore $\frac{2}{\gamma} = \frac{1}{\alpha}
+ \frac{1}{\beta}$, which proves that $\gamma$ is the harmonic
mean of $\alpha$ and $\beta$.  Finally, by going back to the
equation $(\beta - \gamma)(2\beta + \gamma) = 0$ we consider the
other possibility, namely $\beta- \gamma =0$; $\beta = \gamma$
(note that $\beta - \gamma$ and $2\beta + \gamma$ cannot both be
zero for this would imply $\beta = 0$, violating the problem's
assumption that $\beta \neq 0$).  Since $\beta = \gamma$ and
$2\alpha = \beta + \gamma$, we conclude $\alpha = \beta = \gamma$.
And then trivially, $\frac{2}{\gamma} = \frac{1}{\alpha} +
\frac{1}{\beta}$, so we are done.

\item[P15.] We partition the set of natural numbers in disjoint
classes or groups as follows:
$\{1\},\{2,3\},\{4,5,6\},\{7,8,9,10\},\ldots $; the $n$th class
contains $n$ consecutive positive integers starting with
$\frac{n\cdot(n-1)}{2} + 1$.  Find the sum of the members of the
$n$th class.

\vspace{.15in}

\noindent{\bf Solution:}  First let us make clear why the first
member of $n$th class is the number $\frac{n(n-1)}{2} +1$; observe
that the $n$th class is preceded by $(n-1)$ classes; so since the
$k$th class, $1 \leq k \leq n-1$, contains exactly $k$ consecutive
integers, then there precisely $(1+2+\ldots +k+\ldots +(n-1))$
consecutive natural numbers preceding the $n$th class; but the sum
$1+2+\ldots + (n-2)+(n-1)$ is the sum of the first $(n-1)$ terms
of the infinite arithmetic progression that has first term $a_1
=1$ difference $d=1$, hence

$$
\begin{array}{rcl}
1+2+\ldots+(n-1)& = & a_1+a_2+\ldots+a_{n-1} = \frac{(n-1)\cdot(a_1+a_{n-1})}{2}\\
\\
& = & \frac{(n-1)(1+(n-1))}{2} = \frac{(n-1)\cdot n}{2}.\end{array}
$$

\noindent This explains why the $n$th class starts with the
natural number $\frac{n(n-1)}{2}+1$; the members of the $n$th
class are the numbers $\frac{n(n-1)}{2} + 1, \ \frac{n(n-1)}{2} +
2, \ldots, \frac{n(n-1)}{2}+n$.  These $n$ numbers form a finite
arithmetic progression with first term
$\underset{a}{\underbrace{\frac{n(n-1)}{2}+1}}$ and difference
$d=1$.  Hence their sum is equal to

$$
\begin{array}{rcl}
\frac{n\cdot [2a+(n-1)d]}{2}&  =& \frac{n\cdot\left[2\left(\frac{n(n-1)}{2}+1\right)+(n-1)\right]}{2}\\
\\
& = &  \frac{n\cdot [n(n-1)+2+n-1]}{2} = \frac{n\cdot [n^2-n+2+n-1]}{2} = \fbox{$\frac{n\cdot(n^2+1)}{2}$}\end{array}
$$

\item[P16.]  We divide 8,000 objects into $(n+1)$ groups of which
the first $n$ of them contain $5,8,11,14,\ldots ,[5+3\cdot(n-1)]$
objects respectively;  and the $(n+1)$th group contains fewer than
$(5+3n)$ objects; find the value of the natural number $n$ and the
number of objects that the $(n+1)$th group contains.

\vspace{.15in}

\noindent{\bf Solution:}  The total number of objects that first
$n$ groups contain is equal to, $S_n=5+8+11+14+ \ldots +
[5+3(n-1)]$; this sum, $S_n$, is the sum of the first $n$ terms of
the infinite arithmetic progression with first term $a_1 = 5$ and
difference $d = 3$; so that its $n$th term is $a_n=5+3(n-1)$.
According to Theorem 2, $S_n = \frac{n\cdot[a_1 + a_n]}{2} =
\frac{ n\cdot[5+5+3(n-1)]}{2} = \frac{n\cdot [5+5+3n-3]}{2} =
\frac{n\cdot (7+3n)}{2}$.  Thus, the $(n+1)$th group must contain,
$8,000 - \frac{n(7+3n)}{2}$ objects.  By assumption, the $(n+1)$th
group contains fewer than $(5+3n)$ objects.  Also $8,000 -
\frac{n(7+3n)}{2}$ must be a nonnegative integer, since it
represents the number of objects in a set (the $(n+1)$th class;
theoretically this number may be zero).  So we have two
simultaneous inequalities to deal with:

$$
0 \leq 8,000 - \frac{n(7+3n)}{2} \Lra \frac{n(7+3n)}{2} \leq 8,000;\ \ n(7+3n)\leq 16,000.
$$

\noindent And (the other inequality)

$$
\begin{array}{rcl}8,000 - \frac{n(7+3n)}{2} & < & 5 + 3n \Lra 16,000 - n(7+3n)<10+6n \Lra 16,000\\
&  <  &3n^2 + 13n + 10 \Lra 16,000 < (3n+10)(n+1).\end{array}
$$

\noindent So we have the following system of two simultaneous
inequalities

$$
\left.\begin{array}{rc}  & n(7+3n) \leq 16,000\\
\\
{\rm and} & 16,000 < (3n+10)(n+1) \end{array}\right\} \begin{array}{c}(1)\\
\\
(2)\end{array}
$$

\noindent Consider (1): At least one of the factors $n$ and $7+3n$
must be less than or equal to $\sqrt{16,000}$; for if both were
greater than $\sqrt{16,000}$ then their product would exceed
$\sqrt{16,000}\cdot\sqrt{16,000} = 16,000$, contradicting
inequality (1); and since $n<7+3n$, it is now clear that the
natural number $n$ cannot exceed $\sqrt{16,000}:n\leq
\sqrt{16,000}\Lra n \leq \sqrt{16\cdot 10^3};\ n \leq 4\cdot
\sqrt{10^2\cdot 10};\ n \leq 4 \cdot 10\cdot \sqrt{10} =
40\sqrt{10}$ so $40\sqrt{10}$ is a necessary upper bound for $n$.
The closest positive integer to $40\sqrt{10}$, but less than
$40\sqrt{10}$ is the number $126$; but actually, an upper bound
for $n$ must be much less than $126$ in view of the factor $7+3n$.
If we consider (1), we have $3n^2+7n-16,000 \leq 0$ (3)

The two roots of the quadratic equation $3x^2+7x - 16,000 = 0$ are
the real numbers $r_1 = \frac{-7+\sqrt{(7)^2-4(3)(-16,000)}}{6} =
\frac{-7+\sqrt{192,049}}{6} = \approx$\linebreak $ {\rm
approximately}\ 71.872326$; and $r_2 = \frac{-7-\sqrt{192,049}}{6}
\approx -74.20566$.

Now, it is well known from precalculus that if $r_1$ and $r_2$ are
the two roots of the quadratic polynomial $ax^2+bx+c$, then $ax^2
+ bx+c = a\cdot(x-r_1)(x-r_2)$, for all real numbers $x$.  In our
case $3x^2+7x-16,000 = 3\cdot(x-r_1)(x-r_2)$, where $r_1$ and
$r_2$ are the above calculated real numbers.  Thus, in order for
the natural number $n$ to satisfy the inequality (3),
$3n^2+7n-16,000 \leq 0$;  it must satisfy $3(n-r_1)(n-r_2)\leq 0$;
but this will only be true if, and only if, $r_1 \leq n \leq r_2$;
$-74.20566 \leq n \leq 71.872326$; but $n$ is a natural number;
thus $1 \leq n \leq 71$;  this upper bound for $n$ is much lower
than the upper bound of the upper bound $126$ that we estimated
more crudely earlier.  Now consider inequality (2): it must hold
true simultaneously with (1); which means we have,

$$
\left.\begin{array}{rl}
& 16,000 < (3n+10)\cdot(n+1)\\
\\
{\rm and} & 1 \leq n \leq 71 \end{array} \right\}
$$

\noindent If we take the highest value possible for $n$; namely $n
= 71$, we see that $(3n+10)(n+1) = (3\cdot(71)+10)\cdot(72)=
(223)(72) = 16,052$ which exceeds the number $16,000$, as desired.
But, if we take the next smaller value, $n = 70$, we have
$(3n+10)(n+1) = (220)(71) = 15,620$ which falls below $16,000$.
Thus, this problem has a unique solution, \fbox{$n = 71$}.  The
total number of objects in the first $n$ groups (or 71 groups) is
then equal to,

$$
\frac{n\cdot(7+3n)}{2} = \frac{(7)\cdot(7+3(7))}{2} = \frac{(71)\cdot(220)}{2} = (71)\cdot(110) - 7,810.
$$

\noindent Thus, the $(n+1)$th or $72$nd group contains,  $8,000 -
7,810 = \fbox{190}$ objects;  note that $190$ is indeed less that
$5n+3=5(71)+3 = 358$.

\item[P17.]
\begin{enumerate}
\item[(a)]  Show that the real numbers
$\frac{\sqrt{2}+1}{\sqrt{2}-1},\ \frac{1}{2-\sqrt{2}},\
\frac{1}{2}$, can be three consecutive terms of a geometric
progression.  Find the ratio $r$ of any geometric progression that
contains these three numbers as consecutive terms.

\item[(b)] Find the value of the infinite sum of the terms of the
(infinite) geometric progression whose first three terms are the
numbers $\frac{\sqrt{2}+1}{\sqrt{2}-1},\ \frac{1}{2-\sqrt{2}},\
\frac{1}{2}; \ \left(\frac{\sqrt{2}+1}{\sqrt{2}-1} \right) +
\left( \frac{1}{2-\sqrt{2}} \right)+ \frac{1}{2} + \ldots$ .
\end{enumerate}
\vspace{.15in}

\noindent {\bf Solution:}
\begin{enumerate}
\item[(a)]  Apply Theorem 7:  the three numbers will be consecutive terms of a geometric progression if, and only if,

\vspace{.15in}

\hspace{.5in}$\left({\displaystyle \frac{1}{2-\sqrt{2}}}\right)^2 = {\displaystyle \frac{(\sqrt{2}+1)}{(\sqrt{2} - 1)}} \cdot \frac{1}{2}$ \hfill (1)

\vspace{.15in}

\noindent Compute the left-hand side:

$$\begin{array}{rcl}{\displaystyle \frac{1}{(2-\sqrt{2})^2}} & = &{\displaystyle \frac{1}{4-4\sqrt{2}+2}} = {\displaystyle \frac{1}{6-4\sqrt{2}}}\\
\\
& = &{\displaystyle \frac{1}{2(3-2\sqrt{2})}} = {\displaystyle \frac{3+2\sqrt{2}}{2\cdot(3-2\sqrt{2})(3+2\sqrt{2})}}\\
\\
& =& {\displaystyle \frac{3+2\sqrt{2}}{2\cdot[9-8]}} = {\displaystyle \frac{3+2\sqrt{2}}{3}}.\end{array}
$$

\noindent Now we simplify the right-hand side:

$$\begin{array}{rcl}\left( {\displaystyle \frac{\sqrt{2}+1}{\sqrt{2}-1}}\right) \cdot {\displaystyle \frac{1}{2}}
 & = &{\displaystyle \frac{1}{2} \cdot \frac{(\sqrt{2}+1)^2}{(\sqrt{2}-1)(\sqrt{2}+1)}}\\
\\
&  = & {\displaystyle \frac{1}{2} \cdot \frac{(2+2\sqrt{2} + 1)}{(2-1)} = \frac{3+\sqrt{2}}{2}} \end{array}
$$

\noindent so the two sides of (1) are indeed equal; (1) is a true
statement.  Thus, the three numbers can be three consecutive terms
in a geometric progression.  To find $r$, consider
$\left(\frac{\sqrt{2}+1}{\sqrt{2}-1}\right) \cdot r =
\frac{1}{2-\sqrt{2}}$; and also $\left(\frac{1}{2-\sqrt{2}}\right)
\cdot r = \frac{1}{2}$;  from either of these two equations we can
get the value of $r$; if we use the second equation we have,
\fbox{$r = \frac{2-\sqrt{2}}{2}$}. \item[(b)] Since $|r| = \left|
\frac{2-\sqrt{2}}{2} \right| = \frac{2-\sqrt{2}}{2} < 1$,
according to Remark 6, the sum $a + ar + ar^2 + \ldots +
ar^{n-1}+\ldots$ converges to $\frac{a}{1-r}$; in our case $a =
\frac{\sqrt{2}+1}{\sqrt{2}-1}$ and $ r = \frac{2-\sqrt{2}}{2}$.
Thus the value of the infinite sum is equal to

$$\begin{array}{rcl}{\displaystyle \frac{a}{1-r}} & = &{\displaystyle
\frac{\frac{\sqrt{2}+1}{\sqrt{2}-1}}{1-\left(\frac{2-\sqrt{2}}{2}\right)}}=
{\displaystyle \frac{\frac{\sqrt{2}+1}{\sqrt{2}-1}}{\frac{2-\left(2-\sqrt{2}\right)}{2}}}\\
\\
 & = & {\displaystyle \frac{2(\sqrt{2}+1)}{\sqrt{2}(\sqrt{2}-1)} =
 \frac{2(\sqrt{2}+1)\cdot (\sqrt{2}+1)\cdot \sqrt{2}}{\sqrt{2} \cdot \sqrt{2} \cdot (\sqrt{2}-1)(\sqrt{2}+1)}} \\
\\
& =& {\displaystyle  \frac{2\sqrt{2} \cdot (\sqrt{2}+1)^2}{2\cdot (2-1)}=  \sqrt{2} \cdot (2+2\sqrt{2}+1) = \sqrt{2} \cdot (3 + 2\sqrt{2})}\\
\\
& = & 3\sqrt{2} + 2 \cdot \sqrt{2} \cdot \sqrt{2} = 3\sqrt{2} + 4 = \fbox{$4+3\sqrt{2}$}
\end{array}
$$

\item[P18.]  (For student who had Calculus.)  If $|\rho | < 1$ and $|\beta\rho | < 1$, calculate the infinite sum,

$$S = \underset{1{\rm st}}{\underbrace{\alpha\rho}} +
\underset{2{\rm nd}}{(\underbrace{\alpha+\alpha\beta})} \rho^2 + \ldots +
\underset{n{\rm th\ term}}{(\underbrace{\alpha +\alpha\beta + \ldots + \alpha \beta^{n-1}})} \rho^n + \ldots \ .
$$

\noindent{\bf Solution:}  First we calculate the $n$th term which itself is a sum of $n$ terms:

$$
(\alpha+\alpha\beta + \ldots + \alpha \beta^{n-1}) \cdot \rho^n = \alpha \cdot \rho^n \cdot
(1+\beta +\ldots + \beta^{n-1}) = \alpha \cdot \rho^n \cdot \left(\frac{\beta^n-1}{\beta-1}\right)
$$

\noindent by Theorem 5(ii). Now we have,

$$
\begin{array}{rcl}
S & = & \alpha \rho + (\alpha + \alpha \beta)\rho^2 + \ldots + \alpha \cdot \rho^n \cdot
\left( \frac{\beta^n-1}{\beta-1}\right) + \ldots \\
\\
S & = & \alpha \rho \left( \frac{\beta-1}{\beta-1}\right) + \alpha \rho^2 \cdot \left(\frac{\beta^2-1}{\beta-1}\right)\\
\\
& & + \ldots + \alpha \cdot \rho^n \cdot \left( \frac{\beta^n-1}{\beta -1}\right) + \ldots
\end{array}
$$

\noindent Note that $S = {\displaystyle \lim_{n\ra \infty}} S_n$, where

$$
\begin{array}{rcl}
S_n & = & \alpha \rho \cdot \left(\frac{\beta -1}{\beta - 1}\right) + \alpha \rho^n \cdot
\left( \frac{\beta^2 - 1}{\beta - 1}\right) + \ldots + \alpha \cdot \rho^n \cdot \left(\frac{\beta^n-1}{\beta-1}\right);\\
\\
S_n & = & \left(\frac{\alpha\rho}{\beta -1}\right) \left[ (\beta - 1) + \rho (\beta^2 -1) +
\ldots + \rho^{n-1} \cdot (\beta^n - 1)\right]\\
\\
S_n & = & \left(\frac{\alpha\rho}{\beta -1}\right) \left[ \beta \cdot [1 + (\rho \beta) +
\ldots + (\rho \beta)^{n-1}] - (1 + \rho + \ldots + \rho^{n-1})\right]\\
\\
S_n & = & \left( \frac{\alpha\rho}{\beta -1}\right) \cdot \left[ \beta \cdot
\frac{[(\rho\beta)^n-1]}{\rho \beta -1} - \left( \frac{\rho^n-1}{\rho  -1}\right)\right]
\end{array}
$$

\noindent Now, as $n\ra \infty$, in virtue of $|\rho \beta |<1$
and $|\rho | < 1$ we have, ${\displaystyle \lim_{n\ra \infty}}
\frac{[(\rho \beta)^n-1]}{\rho \beta - 1} = \frac{-1}{\rho \beta -
1} = \frac{1}{1 - \rho \beta}$ and ${\displaystyle \lim_{n\ra
\infty}} \frac{\rho^n -1}{\rho - 1} = \frac{1}{1-\rho}$.  Hence,

$$
\begin{array}{rcl}
S & = & {\displaystyle \lim_{n\ra \infty}} S_n = \left( \frac{\alpha \rho}{\beta -1}\right) \cdot
\left[ \beta \cdot \left( \frac{1}{1-\rho\beta}\right)- \left( \frac{1}{1-\rho}\right)\right];\\
\\
S& = & \left( \frac{\alpha \rho}{\beta - 1}\right) \cdot \left[ \frac{\beta (1-\rho)-(1- \rho \beta)}{(1 - \rho \beta)
\cdot (1-\rho)}\right] = \frac{\alpha \rho \cdot (\beta -1 )}{(\beta -1) \cdot (1 -\rho\beta)(1-\rho)};\end{array}
$$

\hspace*{.45in}\fbox{$\begin{array}{c}S = \frac{\alpha \rho}{(1-\rho \beta) \cdot (1-\rho)}\end{array}$}
\end{enumerate}

\item[P19.] Let $m, n$ and $\ell$ be distinct natural numbers; and $a_1 , \ldots , a_k , \ldots$, an
infinite arithmetic progression with first nonzero term $a_1$ and difference $d$.
\begin{enumerate}
\item[(a)]  Find the necessary conditions that $n,\ell$, and $m$ must satisfy in order that,
$$\underset{\underset{{\rm first}\ m\ {\rm terms}}{\rm sum\ of\ the}}{\underbrace{a_1 +a_2 + \ldots + a_m}} =
\underset{\underset{{\rm next}\ n\ {\rm terms}}{\rm sum\ of\ the}}{\underbrace{a_{m+1}+\ldots + a_{m+n}}} =
 \underset{\underset{{\rm next}\ \ell\ {\rm terms}}{\rm sum\ of\ the}}{\underbrace{a_{m+1}+\ldots +a_{m + \ell}}}$$

\item[(b)] If the three sums in part (a) are equal, what must be the relationship between $a_1$ and $d$?
\item[(c)]  Give numerical examples.
\end{enumerate}

\vspace{.15in}

\noindent{\bf Solution:}
\begin{enumerate}
\item[(a)]  We have two simultaneous equations,

$
\left.\begin{array}{rl}
& a_1 + a_2 + \ldots + a_m = a_{m+1} +\ldots +a_{m+n} \\
{\rm and }\\
& a_{m+1}+ \ldots + a_{m+n} = a_{m+1} + \ldots + a_{m+\ell} \end{array} \right\}
$ \hfill (1)

\noindent According to Theorem 2 we have,

$$\begin{array}{rrcl}
&a_1 + a_2 + \ldots + a_m & =&  \frac{m\cdot [2a_1 + (m-1)d]}{2};\\
\\
& a_{m+1}+\ldots + a_{m+n} & = & \frac{n\cdot[a_{m+1}+a_{m+n}]}{2}\\
\\
&& = & \frac{ n\cdot[(a_1 + md) + (a_1 + (m+n-1)d)]}{2}\\
\\
& & = & \frac{n\cdot [2a_1 + (2m+n-1)d]}{2};\\

{\rm and}& & & \\
& a_{m+1}+\ldots + a_{m+\ell} & = & \frac{\ell \cdot[2a_1 + (2m+\ell -1)d]}{2}
\end{array}
$$

\noindent Now let us use the first equation in (1):

$$\begin{array}{rcl}
\frac{m\cdot [2a_1 + (m-1)d]}{2} & = & \frac{n\cdot[2a_1 + (2m+n-1)d]}{2};\\
\\
2ma_1 +m(m-1)d & = & 2na_1 + n\cdot (2m+n-1)d;\\
\\
2a_1 \cdot (m-n) & = & [n \cdot (2m + n-1) - m(m-1)]d;\\
\\
2a_1 \cdot (m-n) & = & [ 2nm + n^2 - m^2 + m -n ]d;
\end{array}
$$

According to hypothesis $a_1 \neq 0$ and $m-n \neq 0$; so the right-hand side must also be nonzero and,

\vspace{.15in}

\hspace*{1.0in}$d = \frac{2a_1 \cdot (m-n)}{2nm+n^2 - m^2 + m-n}$ \hfill (2)

\vspace{.15in}

Now use the second equation in (1):

$$
\begin{array}{rrcl}
& \frac{n\cdot[2a_1 + (2m+n-1)d]}{2} & = & \frac{\ell \cdot [2a_1 + (2m + \ell -1)d]}{2}\\
\\
\Lra & 2na_1 + n(2m+n-1)d & = & 2\ell a_1 + \ell (2m+\ell -1)d\\
\\
\Lra & 2a_1 \cdot (n-\ell) & = & [ \ell (2m+\ell -1) - n(2m+n-1)]d\\
\\
\Lra & 2a_1 \cdot (n-\ell ) & = & [2m \cdot (\ell -n)+(\ell^2 - n^2)-(\ell - n)]d\\
\\
\Lra & 2a_1 \cdot (n-\ell ) & = &[2m \cdot (\ell - n) + (\ell - n)(\ell + n) - (\ell - n)]d\\
\\
\Lra & 2a_1 \cdot (n- \ell ) & = & (\ell - n) \cdot [2m + \ell + n-1]d;
\end{array}
$$

\noindent and since $n - \ell \neq$, we obtain $-2a_1 = (2m+\ell _n-1)d$;

\vspace{.15in}

\hspace*{1.5in}$d = \frac{-2a_1}{2m+\ell + n-1}$ \hfill (3)

\vspace{.15in}

\noindent (Again, in virtue of $a_1 \neq 0$, the product $(2m + \ell + n-1)d$ must also be nonzero,
so $2m + \ell + n-1 \neq 0$, which is true anyway since, obviously, $2m + \ell + n$ is a natural number greater than 1).

Combining Equations (2) and (3) and cancelling out the factor $2a_1 \neq 0$ from both sides we obtain,

$$
\frac{m-n}{2nm+n^2 - m^2 + m-n} = \frac{-1}{2m + \ell + n-1}
$$

\noindent Cross multiplying we now have,

$$
\begin{array}{cl}
&(m-n)\cdot (2m + \ell + n-1)\\
\\
 = & (-1) \cdot (2nm + n^2 - m^2 + m - n);\\
\\
& 2m^2 + m\ell + mn -m - 2mn - n\ell - n^2 + n\\
\\
 = & -2mn - n^2 + m^2 - m + n;\\
\\
& m^2 + m\ell - n\ell + mn  =  0.
\end{array}
$$
\noindent We can solve for $n$ in terms of $m$ and $\ell$ (or for $\ell$ in terms of $m$ and $n$) we have,
$$
n\cdot (\ell - m) = m\cdot (m+\ell ) \Ra \fbox{$n = \frac{m\cdot(m + \ell)}{\ell - m}$},\ {\rm since}\  \ell - m \neq 0.
$$

Also, we must have \fbox{$\ell > m$}, in view of the fact that $n$
is a natural number and hence positive (also note that these two
conditions easily imply $n > m$ as well).  But, there is more: The
natural number $\ell - m$ must be a divisor of the product $m
\cdot (m+\ell)$.  Thus, the conditions are:

\begin{enumerate}
\item[(A)] $\ell > m$
\item[(B)] $(\ell - m)$ is a divisor of $m \cdot (m+\ell )$ and
\item[(C)] $n = \frac{m\cdot (m+\ell)}{\ell - m}$
\end{enumerate}

\item[(b)] As we have already seen $d$ and $a_1$ must satisfy both
conditions (2) and (3).  However, under conditions (A), (B), and
(C), the two conditions (2) and (3) are, in fact, equivalent, as
we have already seen;  so $d = \frac{-2a_1}{2m+\ell + n-1}$
(condition (3)) will suffice.

\item[(c)] Note that in condition (C), if we choose $m$ and $\ell$
such $\ell - m$ is \ positive and $(\ell - m)$ is a divisor of
$m$, then clearly the number $n = \frac{m\cdot (m+\ell)}{\ell -
m}$, will be a natural number.  If we set $\ell - m = t$, then
$m+\ell = t + 2m$, so that

$$
n= \frac{m\cdot(t + 2m)}{t} = m + \frac{2m^2}{t}.
$$

\noindent So if we take $t$ to be a divisor of $m$, this will be
sufficient for $\frac{2m^2}{t}$ to be a positive integer.  Indeed,
set $m = M\cdot t$, then $n = M\cdot t + \frac{2M^2t^2}{t} =
M\cdot t + 2M^2 \cdot t = t\cdot M \cdot (1 + 2M)$.  Also, in
condition (3) , if we set $a_1 = a$, then (since $\ell = m+t =
Mt+t$)

$
\begin{array}{rcl}
d & = & \frac{-2a}{2M\cdot t+(Mt+t)+Mt + 2M^2 t - 1};\\
\\
d & = & \frac{-2a}{4Mt + t + 2M^2 t -1} .\end{array}$ \hfill (4)

\noindent Thus, the formulas $\ell = Mt+t,\ n = Mt+2M^2 \cdot t$
and (4) will generate, for each pair of values of the natural
numbers $M$ and $t$, an arithmetic progression that satisfies the
conditions of the problem; for any nonzero value of the first term
$a$.

\noindent{\bf Numerical Example:}  If we take $t=3$ and $M=4$, we then have $m = M\cdot t = 3 \cdot 4 = 12;\ n =
 t\cdot M \cdot (1+2M) = 12 \cdot (1 + 8) = 108$, and $\ell = m + t = 12+3=15$.  And,

$$d= \frac{-2a}{2m+\ell+n-1} = \frac{-2a}{24+15 + 108-1} = \frac{-2a}{146} = \frac{-a}{73}.
$$

Now let us compute
$$\begin{array}{rcl}a_1 + \ldots + a_m& =& {\displaystyle \frac{m\cdot [2a + (m-1)d]}{2} = \frac{12\cdot
\left[2a+11 \cdot \left(\frac{-a}{73}\right)\right]}{2}}\\
\\
& =& {\displaystyle \frac{12\cdot [146a - 11a]}{2\cdot 73} = \frac{6\cdot(135a)}{73}} =
 {\displaystyle \frac{810a}{73}}.\end{array}
$$
\noindent Next,
$$\begin{array}{rl}
 & a_{m+1}+\ldots + a_{m+n^{\prime}}\\
\\
 =&  \frac{n\cdot [2a+(m+n-1)d]}{2}\\
\\
 = & \frac{108 \cdot \left[2a+(24+108-1)\cdot \left(\frac{-a}{73}\right)\right]}{2}\\
\\
= & \frac{108}{2} \cdot \frac{[146a -131 a]}{73}\\
\\
= &\frac{(54)(15a)}{73} = \frac{810a}{73}\end{array}
$$
\noindent and
$$ \begin{array}{cl}
& a_{m+1} + \ldots + a_{m+\ell} \\
\\
= & {\displaystyle \frac{ \ell \cdot[2a + (2m+\ell-1)d]}{2}}\\
\\
= &  {\displaystyle\frac{15\cdot\left[2a+(24+15-1) \cdot \left(\frac{-a}{73}\right)\right]}{2}} \\
\\
= & {\displaystyle \frac{15}{2} \cdot \frac{[146a - 38a]}{73}}\\
\\
= & {\displaystyle \frac{15}{2} \cdot \frac{(108)a}{73} = \frac{(15)(54a)}{73}}\\
\\
= & {\displaystyle \frac{810a}{73}}.
\end{array}
$$
\noindent Thus, all three sums are equal to $\frac{810a}{73}$.
\end{enumerate}

\item[P20.]  If the real numbers $a,b,c$ are consecutive terms of an arithmetic progression and
$a^2, b^2, c^2$ are consecutive terms of a harmonic progression, what conditions must the numbers
$a,b,c$ satisfy?  Describe all such numbers $a,b,c$.

\vspace{.15in}

\noindent{\bf Solution:}  By hypothesis, we have
$$2b=a + c \ {\rm and}\ \frac{2}{b^2} = \frac{1}{a^2} + \frac{1}{c^2}
$$
\noindent so $a,b,c$ must all be nonzero real numbers.  The second equation is equivalent
to $b^2 = \frac{2a^2c^2}{a^2+c^2}$ and $abc \neq 0$; so that, $b^2 (a^2+c^2) = 2a^2c^2 \Lra b^2 \cdot [(a+c)^2 - 2ac] = 2a^2 c^2$.
 Now substitute for $a+c = 2b$:

$$\begin{array}{rl}
& b^2 \cdot [(2b)^2 - 2ac] = 2a^2c^2 \\
\\
\Lra & 4b^4 - 2acb^2 - 2a^2c^2 = 0;\\
\\
& 2b^4 - acb^2 - a^2 c^2 = 0
\end{array}
$$

\noindent At this stage we could apply the quadratic formula since
$b^2$ is a root to the  equation $2x^2 -acx-a^2c^2 = 0$; but the
above equation can actually be factored. Indeed,

\vspace{.15in}

\hspace*{1.0in}$\begin{array}{rcl}
b^4 - acb^2 + b^4 - a^2c^2 &= & 0;\\
\\
b^2(b^2 - ac) + (b^2)^2 - (ac)^2 & = & 0;
\end{array}
$

\hspace*{.52in}$\begin{array}{rcl}
b^2 \cdot (b^2 - ac) + (b^2 - ac)(b^2 + ac) & = & 0;\\
(b^2 - ac) \cdot (2b^2 + ac) & = & 0 \end{array}
$\hfill (1)

\vspace{.15in}

According to Equation (1), we must have $b^2 - ac= 0$; or
alternatively $2b^2 + ac = 0$.  Consider the first possibility,
$b^2 - ac = 0$.  Then, by going back to equation $\frac{2}{b^2}=
\frac{1}{a^2} + \frac{1}{c^2}$ we obtain $\frac{2}{ac} =
\frac{1}{a^2} + \frac{1}{c^2} \Lra \frac{2a^2c^2}{ac} = a^2 + c^2
\Lra 2ac = a^2 + c^2$; $a^2 + c^2 - 2ac = 0 \Lra (a-c)^2 = 0$; $a
= c$ and thus $2b = a+c$ implies $b = a = c$.

Next, consider the second possibility in Equation (1):  $2b^2 + ac
= 0 \Lra 2b^2 = -ac$;  which clearly imply that one of $a$ and $c$
must be positive, the other negative.  Once more going back to

\vspace{.15in}

\hspace*{.75in}$\begin{array}{rl}
& \frac{2}{b^2} = \frac{1}{a^2} + \frac{1}{c^2};\ \frac{4}{2b^2} = \frac{1}{a^2} + \frac{1}{c^2}\\
\\
 \Lra & \frac{4}{-ac} = \frac{c^2+a^2}{a^2c^2 }\\
\\
 \Lra & -4ac = c^2 + a^2$; $a^2 + 4ac + c^2 = 0
\end{array}$ \hfill (2)

\vspace{.15in}

Let $t = \frac{a}{c};\ a=c\cdot t$ then Equation (2) yields (since $ac \neq 0$),

\vspace{.15in}

\hspace*{1.5in}$t^2 + 4t+1 = 0$ \hfill (3)

\vspace{.15in}

\noindent Applying the quadratic formula to Equation (3), we now have

$$
\begin{array}{l}t = {\displaystyle \frac{-4\pm \sqrt{16-4}}{2} = \frac{-4 \pm 2\sqrt{3}}{2}};\\
\\
t = -2 \pm \sqrt{3};
\end{array}
$$

\noindent note that both numbers $-2+\sqrt{3}$ and $-2-\sqrt{3}$
are negative and hence both acceptable as solutions, since we know
that $a$ and $c$ have opposite sign, which means that $t =
\frac{a}{c}$ must be negative.  So we must have either $a = (-2 +
\sqrt{3})c$; or alternatively $a = -(2+\sqrt{3})\cdot c$.  Now, we
find $b$ in terms of $c$.  From $2b^2 = -ac$; $b^2 =
-\frac{ac}{2}$; note that the last equation says that either the
numbers $- \frac{a}{2}, b, c$ are the successive terms of a
geometric progression; or the numbers $-a,b,\frac{c}{2}$ (or any
of the other two possible permutations: $a,b,-\frac{c}{2}$,
$\frac{a}{2}, b, -c$; and four more that are obtained by switching
$a$ with $c$).  So, if $a = (-2+\sqrt{3})c$, then from $2b = a +
c;\ b = \frac{a+c}{2} = \frac{(-2+\sqrt{3})c+c}{2} =
\frac{(\sqrt{3}-1)c}{2}$.  And if $a = -(2+\sqrt{3})c,\ b =
\frac{a+c}{2} = \frac{-(2+\sqrt{3}) c+c}{2} =
\frac{-(1+\sqrt{3})c}{2}$.  So, in conclusion we summarize as
follows:

Any three real numbers $a,b,c$ such that $a,b,c$ are consecutive
terms of an arithmetic progression  and $a^2,b^2,c^2$ the
successive terms of a harmonic progression must fall in exactly
one of three classes:

\begin{enumerate}
\item[(1)] $a = b = c;\ c$ can be any nonzero real number
\item[(2)] $a = (-2+\sqrt{3})\cdot c,\ b =
\frac{(\sqrt{3}-1)c}{2};\ c$ can be any positive real;

\item[(3)] $a = (2 + \sqrt{3})c,\ b= \frac{-(1+\sqrt{3})}{2}c;\ c$
can be any positive real.
\end{enumerate}

\item[P21.] Prove that if the positive real numbers $\alpha,
\beta, \gamma$ are consecutive  members of a geometric
progression, then $\alpha ^k + \gamma ^k \geq 2 \beta^k$, for
every natural number $k$.

\vspace{.15in}

\noindent{\bf Solution:}  Given any natural number $k$, we can
apply the arithmetic-geometric  mean inequality of Theorem 10,
with $n = 2$, and $a_1 = \alpha^k,\ a_2 = \gamma^k$, in the
notation of that theorem:

$$
\frac{\alpha^k + \gamma^k}{2} \geq \sqrt{\alpha^k \cdot \gamma^k} = \sqrt{(\alpha\gamma)^k}.
$$

\noindent But since $\alpha ,\beta , \gamma$ are consecutive terms
of a geometric progression, we must also have $\beta^2 = \alpha
\gamma$.  Thus the above inequality implies,

$$
\begin{array}{rccl}
& \frac{\alpha^k + \gamma^k}{2} & \geq & \sqrt{(\beta^2)^k};\\
& \frac{\alpha^k + \gamma^k}{2} & \geq & \sqrt{(\beta^k)^2} \\

\Ra & \frac{\alpha^k + \gamma^k}{2} & \geq & \beta ^k\\
\Ra & \alpha^k + \gamma^k & \geq & 2\beta^k,
\end{array}
$$

\noindent and the proof is complete.
\end{enumerate}

\section{\bf Unsolved problems}

\begin{enumerate}
\item[1.] Show that if the sequence $a_1,a_2,\ldots , a_n,\ldots$
, is an arithmetic  progression, so is the sequence $c \cdot a_1,
c\cdot a_2, \ldots, c \cdot a_n, \ldots$ , where $c$ is a
constant. \item[2.]  Determine the difference of each arithmetic
progression which has first term $a_1 = 6$ and contains the
numbers $62$ and $104$ as its terms. \item[3.] Show that the
irrational numbers $\sqrt{2}, \ \sqrt{3},\ \sqrt{5}$ cannot be
terms of an arithmetic progression. \item[4.] If $a_1, a_2, \ldots
, a_n , \ldots$ is an arithmetic progression and $a_k = \alpha,\
a_m=\beta,\ a_{\ell}= \gamma$, show that the natural numbers
$k,m,\ell$ and the real numbers $\alpha, \beta, \gamma$, must
satisfy the condition

$$
\alpha \cdot (m-\ell) + \beta \cdot (\ell - k) + \gamma \cdot (k - m)=0.
$$

\noindent{\bf Hint:}  Use the usual formula $a_n = a_1 + (n-1)d$,
for $n = k,m,\ell$, to obtain three equations; subtract the first
two and then the last two (or the first and the third) to
eliminate $a_1$; then eliminate the difference $d$ (or solve for
$d$ in each of the resulting equations).

\item[5.]  If the numbers $\alpha,\beta, \gamma$ are successive
terms of an arithmetic progression, then the same holds true for
the numbers $\alpha^2 \cdot(\beta +\gamma),\ \beta^2 \cdot (\gamma
+ \alpha),\ \gamma^2 \cdot (\alpha +\beta)$.

\item[6.]  If $S_k$ denotes the sum of the first $k$ terms of the
arithmetic progression with first term $k$ and difference $d =
2k-1$, find the sum $S_1 + S_2 + \ldots + S_k$.

\item[7.]  We divide the odd natural numbers into groups or
classes as follows:  $\{1\},\{3,5\},\{7,9,11\}, \ldots $ ; the
$n$th group contains $n$ odd numbers starting with $(n \cdot
(n-1)+1)$ (verify this).  Find the sum of the members of the $n$th
group.

\item[8.]  We divide the even natural numbers into groups as
follows:  $\{2\},\{4,6\},$ \linebreak $\{8,10,12\},\ldots$ ; the
$n$th group contains $n$ even numbers starting with $(n(n-1)+2)$.
Find the 'sum of the members of the $n$th group.

\item[9.]  Let $n_1,n_2,\ldots , n_k$ be $k$ natural numbers such
that $n_1 < n_2 < \ldots < n_k$; if the real numbers, $a_{n_1},
a_{n_2},\ldots , a_{n_k}$, are members of an arithmetic
progression (so that the number $a_{n_i}$ is precisely the $n_i$th
term in the progression; $i = 1,2,\ldots ,k)$, show that the real
numbers:

$$
\frac{a_{n_k}-a_{n_1}}{a_{n_2}-a_{n_1}},\ \frac{a_{n_k}-a_{n_2}}{a_{n_2}-a_{n_1}}\  , \ldots ,
\frac{a_{n_k}-a_{n_{k-1}}}{a_{n_2}- a_{n_1}},
$$

\noindent are all rational numbers.

\item[10.]  Let $m$ and $n$ be natural numbers.  If in an
arithmetic progression $a_1,a_2,\ldots , a_k,\ldots $; the term
$a_m$ is equal to $\frac{1}{n}$; $a_m=\frac{1}{n}$, and the term
$a_n$ is equal to $\frac{1}{m};\ a_n = \frac{1}{m}$, prove the
following three statements.

\begin{enumerate}
\item[(a)] The first term $a_1$ is equal to the difference $d$.
\item[(b)] If $t$ is any natural number, then $a_{t\cdot(mn)} =
t$; in other words, the terms $a_{mn},a_{2mn},a_{3mn},\ldots $ ,
are respectively equal to the natural numbers $1,2,3,\ldots$ .
\item[(c)]  If $S_{t\cdot (mn)}$ ($t$ a natural number) denote the
sum of the first $(t\cdot m \cdot n)$ terms of the arithmetic
progression, then $S_{t\cdot(mn)} = \frac{1}{2} \cdot (mn+1)\cdot
t$. In other words, $S_{mn} = \frac{1}{2}(mn+1)$,  $S_{2mn} =
\frac{1}{2} \cdot  (mn+1)\cdot 2,\ S_{3mn} = \frac{1}{2} \cdot
(mn+1)\cdot 3,\ldots $ .
\end{enumerate}

\item[11.]  If the distinct real numbers $a,b,c$ are consecutive terms of a harmonic progression show that

\begin{enumerate}
\item[(a)] $\frac{2}{b} = \frac{1}{b-a} + \frac{1}{b-c}$ and
\item[(b)] $\frac{b+a}{b-a} + \frac{b+c}{b-c} = 2$
\end{enumerate}

\item[12.] If the distinct reals $\alpha, \beta, \gamma $ are
consecutive terms of a harmonic progressionthen the same is true
for the numbers $\alpha, \alpha - \gamma, \alpha - \beta$.

\item[13.]  Let $a = a_1, a_2, a_3, \ldots , a_n, \ldots$ , be a
geometric progression and $k, \ell, m$natural numbers.  If $a_k =
\beta, \ a_{\ell}=\gamma ,\ a_m=\delta$, show that $\beta^{\ell
-m} \cdot \gamma^{m-k} \cdot \delta^{k-\ell} = 1$.

\item[14.]  Suppose that $n$ and $k$ are natural numbers such that
$n > k + 1$; and $a_1 = a, a_2, \ldots , a_t , \ldots $ a
geometric progression, with positive ratio $r \neq 1$, and
positivefirst term $a$.  If $A$ is the value of the sum of the
first $k$ terms of the progression and $B$ is the value of the
last $k$ terms among the $n$ first terms, express the ratio $r$ in
terms of $A$ and $B$ only; and also the first term $a$ in terms of
$A$ and $B$.

\item[15.]  Find the sum $\left( a- \frac{1}{a}\right)^2 +
\left(a^2 - \frac{1}{a^2}\right)^2 + \ldots  \left( a^n -
\frac{a}{a^n}\right)^2$.

\item[16.]  Find the infinite sum $\left( \frac{1}{3} +
\frac{1}{3^2} + \frac{1}{3^3} + \ldots \right)+ \left( \frac{1}{5}
+ \frac{1}{5^2} + \frac{1}{5^3} + \ldots \right)$\linebreak $+
\left( \frac{1}{9} + \frac{1}{9^2} + \frac{1}{9^3} + \ldots
\right) + \ldots +  \underset{k{\rm th\ sum}}{\left(\underbrace{
\frac{1}{(2k+1)} + \frac{1}{(2k+1)^2} + \frac{1}{(2k+1)^3} +
\ldots }\right)} + \ldots $ .

\item[17.] Find the infinite sum $\frac{2}{7} + \frac{4}{7^2} +
\frac{2}{7^3} + \frac{4}{7^4} + \frac{2}{7^5} + \frac{4}{7^6} +
\ldots $ .

\item[18.]  If the numbers $\alpha, \beta, \gamma$ are consecutive
terms of an arithmetic progression and the nonzero numbers $\beta
, \gamma, \delta$ are consecutive terms of a harmonic progression,
show that $\frac{\alpha}{\beta} = \frac{\gamma}{\delta}$.

\item[19.]  Suppose that the positive reals $\alpha, \beta ,\gamma
$ are successive terms of an arithmetic progression and let $x$ be
the geometric mean of $\alpha $ and $\beta$; and let $y$ be the
geometric mean of $\beta$ and $\gamma$.  Prove that $x^2,\beta^2,
y^2$ are successive terms of an arithmetic progression.  Give two
numerical examples.

\item[20.]  Show that if the nonzero real numbers $a,b,c$ are
consecutive terms of a harmonic progression, then the numbers
$a-\frac{b}{2},\ \frac{b}{2},\ c - \frac{b}{2}$, must be
consecutive terms of a geometric progression.  Give two numerical
examples.

\item[21.]  Compute the following sums:

\begin{enumerate}
\item[(i)]  $\frac{1}{2} + \frac{2}{2^2} + \ldots + \frac{n}{2^n}$
\item[(ii)] $1 + \frac{3}{2} + \frac{5}{4} + \ldots + \frac{2n-1}{2^{n-1}}$
\end{enumerate}

\item[22.]  Suppose that the sequence $a_1, a_2, \ldots , a_n,
\ldots$ satisfies  $a_{n+1} = (a_n + \lambda)\cdot \omega$, where
$\lambda $ and $\omega$ are fixed real numbers with $\omega \neq
1$.
\begin{enumerate}
\item[(i)]  Use mathematical induction to prove that for every
natural number, $a_n = a_1 \cdot \omega^{n-1} + \lambda \cdot
\left( \frac{\omega^n - \omega}{\omega -1 }\right)$. \item[(ii)]
Use your answer in part (i) to show that,

$$
\begin{array}{rcl}S_n & = & a_1 + a_2 + \ldots + a_n\\
& =& a_1 \cdot  {\left(\displaystyle \frac{\omega^n - 1}{\omega - 1}\right) + \lambda \cdot
\left( \frac{\omega^{n+1} - n \cdot \omega ^2 + (n-1) \omega}{(\omega - 1)^2}\right)}.\end{array}
$$

\noindent($*$) Such a sequence is called a {\bf semi-mixed} progression.
\end{enumerate}

\item[23.]  Prove part (ii) of Theorem 4.

\item[24.]  Work out part (viii) of Remark 5.

\item[25.]  Prove the analogue of Theorem 4 for geometric
progressions:  if the $(n - m+1)$  positive real numbers
$a_m,a_{m+1}, \ldots , a_{n-1},a_n$ are successive terms of a
geometric progression, then

\begin{enumerate}
\item[(i)]  If the natural number $(n-m+1)$ is odd, then the
geometric mean of the $(n-m+1)$ terms is simply the middle number
$a_{(\frac{m+n}{2})}$.

\item[(ii)]  If the natural number $(n-m+1)$ is even, then the
geometric mean of the $(n-m+1)$ terms must be the geometric mean
of the two middle terms $a_{(\frac{n+m-1}{2})}$ and
$a_{(\frac{n+m+1}{2})}$.
\end{enumerate}
\end{enumerate}

 \end{document}